\documentclass[11pt, amsfonts,amssymb, amsmath]{amsart}

\usepackage{times}
\usepackage{amsfonts}
\usepackage{amssymb} 
\theoremstyle{plain}
\newtheorem{theorem}{Theorem}[section]
\newtheorem{corollary}[theorem]{Corollary}
\newtheorem{lemma}[theorem]{Lemma}
\newtheorem{proposition}[theorem]{Proposition}
\newtheorem{fact}[theorem]{Fact}
\newtheorem{claim}[theorem]{Claim}

\theoremstyle{definition}
\newtheorem{definition}[theorem]{Definition}
\newtheorem{remark}[theorem]{Remark}
\newtheorem{hypothesis}[theorem]{Hypothesis}

\newtheorem{question}[theorem]{Question}

\newtheorem{notation}[theorem]{Notation}
\newtheorem*{theo}{Theorem}

\newcommand{\tp}{\operatorname{tp}}

\newcommand{\K}{\operatorname{\mathcal{K}}}

\newcommand{\cl}{\operatorname{cl}}
\newcommand{\bcl}{\operatorname{bcl}}
\newcommand{\ecl}{\operatorname{ecl}}
\newcommand{\Lstp}{\operatorname{Lstp}}

\newcommand{\PGL}{\operatorname{PGL}}
\newcommand{\ran}{\operatorname{ran}}
\newcommand{\aut}{\operatorname{Aut}}
\newcommand{\Saut}{\operatorname{Saut}}

\newcommand{\infinity}{\infty}

\newcommand{\bpf}{\begin{proof}}
\newcommand{\epf}{\end{proof}}
\renewcommand{\empty}{\emptyset}

\newbox\noforkbox \newdimen\forklinewidth
\setbox0\hbox{$\textstyle\smile$}
\setbox1\hbox to \wd0{\hfil\vrule width \forklinewidth depth-2pt
 height 10pt \hfil}
\setbox\noforkbox\hbox{\lower 2pt\box1\lower 2pt\box0\relax}
\def\unionstick{\mathop{\copy\noforkbox}\limits}

\def\nonfork_#1{\unionstick_{\textstyle #1}}

\setbox0\hbox{$\textstyle\smile$}
\setbox1\hbox to \wd0{\hfil{\sl /\/}\hfil}
\setbox2\hbox to \wd0{\hfil\vrule height 10pt depth -2pt width
               \forklinewidth\hfil}
\newbox\doesforkbox
\setbox\doesforkbox\hbox{\lower 2pt\box1 \lower 2pt\box2\lower2pt\box0\relax}
\def\nunionstick{\mathop{\copy\doesforkbox}\limits}

\def\fork_#1{\nunionstick_{\textstyle #1}}

\begin{document}

\title{Interpreting groups and fields in some nonelementary classes}

\author{Tapani Hyttinen}
\address{Department of Mathematics\\
University of Helsinki\\
P.O. Box 4, 00014 Finland}
\email{tapani.hyttinen@helsinki.fi}

\author{Olivier Lessmann} 
\address{Mathematical Institute\\
Oxford University\\
Oxford, OX1 3LB, United Kingdom}
\email{lessmann@maths.ox.ac.uk}

\author{Saharon Shelah}
\address{Department of Mathematics\\
Rutgers University\\
New Brunswick, New Jersey, United States}
\address{Institute of Mathematics\\
The Hebrew University of Jerusalem\\
Jerusalem 91904, Israel}
\email{shelah@math.huji.ac.il}

\thanks{The first author is partially supported by the Academy of Finland,
grant 40734.
The third author is supported by The Israel Science Foundation,
this is publication 821 on his publication list.}
\date{\today}

\begin{abstract}
This paper is concerned with extensions of geometric stability theory 
to some nonelementary classes.
We prove the following theorem:
\begin{theo}
Let $\mathfrak{C}$ be a large homogeneous model of a stable diagram $D$.
Let $p, q \in S_D(A)$, where $p$ is quasiminimal and $q$ unbounded.
Let $P = p(\mathfrak{C})$ and $Q = q(\mathfrak{C})$.
Suppose that there exists an integer $n < \omega$ such that
\[
\dim(a_1 \dots a_{n}/A \cup C) = n,
\]
for any independent $a_1, \dots, a_{n} \in P$ and finite
subset $C \subseteq Q$,
but
\[
\dim(a_1 \dots a_n a_{n+1} /A \cup C) \leq n,
\]
for some independent $a_1, \dots, a_n, a_{n+1} \in P$ 
and some finite subset $C \subseteq Q$.

Then $\mathfrak{C}$ interprets a group $G$ which acts on the geometry
$P'$ obtained from $P$.
Furthermore, either $\mathfrak{C}$ interprets a non-classical
group, or $n = 1,2,3$ and
\begin{itemize}
\item
If $n = 1$ then $G$ is abelian and acts
regularly on $P'$.
\item
If $n = 2$ the action of $G$ on $P'$ is isomorphic to the affine
action of $K \rtimes K^*$ on the algebraically closed field $K$.
\item
If $n = 3$ the action of $G$ on $P'$ is isomorphic to the action
of $\PGL_2(K)$ on the projective line $\mathbb{P}^1(K)$ of the
algebraically closed field $K$.
\end{itemize}
\end{theo}

We prove a similar result for excellent classes.
\end{abstract}

\maketitle

\section{Introduction}

The fundamental theorem of projective geometry is a striking example
of interplay between geometric and algebraic data:
Let $k$ and $\ell$ distinct lines of, say,
the complex projective plane $\mathbb{P}^2(\mathbb{C})$, with $\infinity$
their point of intersection.
Choose two distinct points $0$
and $1$ on $k \setminus \{ \infinity \}$.
We have the {\em Desarguesian property}: For any 2 pairs 
of distinct points ($P_1, P_2)$ and $(Q_1, Q_2)$
on $k \setminus  \{ \infinity \}$, 
there is an automorphism $\sigma$ of $\mathbb{P}^2(\mathbb{C})$ fixing $\ell$
pointwise, preserving $k$, such that $\sigma(P_i) = Q_i$, for $i = 1,2$.
But for some triples $(P_1, P_2, P_3)$ and $(Q_1, Q_2, Q_3)$ on 
$k \setminus \{ \infinity \}$, this property fails. 
\relax From this, it is possible to endow $k$ with the structure of a division
ring, and another geometric property garantees that it is a field.
Model-theoretically, 
in the language of {\em points} (written $P, Q, \dots$), 
{\em lines} (written $\ell, k, \dots$), and
an {\em incidence relation} $\in$,
we have  
a saturated structure 
$\mathbb{P}^2(\mathbb{C})$, and 
two strongly minimal types $p(x) = \{ x \in k \}$
and $q(x) = \{ x \in \ell \}$. 
The Desarguesian property is equivalent to the following
statement in {\em orthogonality calculus}, which
is the area of model theory dealing with the independent relationship
between types:
$p^2$ is weakly orthogonal to $q^\omega$,
but $p^3$ is not almost orthogonal to $q^\omega$
(see the abstract gives another equivalent condition in terms
of dimension).
\relax From this, we can define a division ring on $k$.
Model theory then gives us more: 
strong minimality guarantees that it is an algebraically
closed field,  
and further conditions that it has characteristic $0$;
it follows that it must be $\mathbb{C}$.

A central theorem of geometric stability, due to Hrushovski~\cite{Hr}
(extending Zilber~\cite{Zi1}),
is a generalisation of this result
to the context of stable first order theory:
Let $\mathfrak{C}$ be a large saturated model of a stable first order theory.
Let $p, q \in S(A)$ be stationary and regular such that for some
$n < \omega$ the type $p^n$ is weakly orthogonal to $q^\omega$ but
$p^{n+1}$ is not almost orthogonal to $q^\omega$.
Then $n=1,2,3$ and if $n=1$ then $\mathfrak{C}$ interprets an abelian
group and if $n =2,3$ then $\mathfrak{C}$ interprets an algebraically closed
field. He further obtains a description of the action for $n=1,2,3$
(see the abstract).  

Geometric stability theory is a branch of first order model theory
that grew out 
of Shelah's classification theory~\cite{Sh:a};
it began with the discovery by Zilber
and Hrushovski that certain model-theoretic problems
(finite axiomatisability of totally categorical 
first order theories~\cite{Zi1},
existence of strictly stable unidimensionaly first order theories~\cite{Hr3})
imposed abstract (geometric) model-theoretic conditions implying the existence
of definable classical groups.  The structure of these groups
was then invoked to solve the problems.
Geometric stability theory
has now developed into a sophisticated body
of techniques which have found remarkable applications
both within model theory (see \cite{Pi}
and \cite{Bu}) and in other areas of mathematics
(see for example the surveys~\cite{Hr1} and \cite{Hr2}).
However, its
applicability is limited at present
to mathematical contexts which are first order axiomatisable.
In order to extend the scope of these techniques, 
it is necessary to develop geometric
stability theory beyond first order logic.
In this paper,
we generalise Hrushovski's result to two non first order settings:
homogeneous model theory
and excellent classes.

{\em Homogeneous model theory}
was initiated by Shelah~\cite{Sh:3}, it consists of studying 
the class of elementary submodels of a large homogeneous,
rather than saturated, model.  Homogeneous model theory is very well-behaved,
with a good notion of stability~\cite{Sh:3}, \cite{Sh:54}, \cite{Hy:4},
\cite{GrLe}, superstability~\cite{HySh}, \cite{HyLe}, 
$\omega$-stability~\cite{Le:1}, \cite{Le:2},
and even simplicity~\cite{BuLe}.  
Its scope of applicability is very broad, as many
natural model-theoretic constructions fit within its framework: 
first order, Robinson theories, existentially closed models, 
Banach space model theory, many
generic constructions,
classes of models with set-amalgamation ($L^n$, infinitary), 
as well as many classical non-first order mathematical
objects like free groups or Hilbert spaces. 
We will consider the stable case
(but note that this context may be unstable from a first
order standpoint),
{\em without} assuming simplicity,
{\em i.e.} without assuming that there is a dependence relation
with all the properties of forking in the first order stable case.
(This contrasts with the work of Berenstein~\cite{Be}, who carries out
some group constructions under the assumption of stability,
simplicity, and the existence of canonical bases.)

{\em Excellence} is a property discovered by 
Shelah~\cite{Sh:87a} and \cite{Sh:87b} 
in his work on categoricity for nonelementary classes:
For example, he proved that, 
under GCH, a sentence in $L_{\omega_1, \omega}$ which is categorical
in all uncountable cardinals is excellent.
On the other hand, excellence 
is central in the classification of almost-free algebras~\cite{MeSh}
and also arises naturally in Zilber's work around complex 
exponentiation~\cite{Zi:2} and \cite{Zi:3}
(the structure $(\mathbb{C}, \exp)$ has intractable 
first order theory since it interprets the integers, but
is manageable in an infinitary
sense). 
Excellence is a condition on the existence of prime models over certain
countable sets (under an $\omega$-stability assumption).
Classification theory for excellent classes is 
quite developed; we have a good understanding
of categoricity (\cite{Sh:87a}, \cite{Sh:87b}, and \cite{Le:3}
for a Baldwin-Lachlan proof), and Grossberg and Hart proved the 
Main Gap~\cite{GrHa}.
Excellence follows   
from uncountable categoricity in the context of homogeneous model theory.
However, excellence is at present restricted to $\omega$-stability
(see \cite{Sh:87a} for the definition),
so excellent classes and stable homogeneous model theory, 
though related, are not comparable.

In both contexts, we lose compactness and saturation, which leads 
us to use various forms of homogeneity instead (model-homogeneity and only
$\omega$-homogeneity in the case
of excellent classes).
Forking is replaced by the appropriate dependence relation, keeping in mind
that not all properties of forking hold at this level
(for example extension and symmetry may fail over certain sets).  
Finally, we have to do without canonical bases.

Each context comes with a notion of monster model $\mathfrak{C}$
(homogeneous or full), which functions as a universal domain;
all relevant realisable types are realised in $\mathfrak{C}$,
and models may be assumed to be submodels of $\mathfrak{C}$.
We consider a {\em quasiminimal} type $p$, {\em i.e.} every definable subset
of its set of realisations in $\mathfrak{C}$
is either bounded or has bounded complement.
Quasiminimal types are a generalisation of strongly minimal types
in the first order case, and play a similar role, for example in 
Baldwin-Lachlan theorems.
We introduce the natural
closure operator on the subsets of $\mathfrak{C}$;
it induces
a pregeometry and a notion of dimension $\dim(\cdot/C)$ 
on the set of realisations of $p$, for any $C \subseteq \mathfrak{C}$.  
We prove:

\begin{theorem}~\label{t:intro1}  Let $\mathfrak{C}$ be a large
homogeneous stable model
or a large full model in the excellent case.
Let $p, q$ be complete types over a finite set $A$,
with $p$ quasiminimal.
Assume that there exists $n < \omega$ such that
\begin{enumerate}
\item
For any independent sequence $(a_0, \dots, a_{n-1})$
of realisations of $p$
and any countable set $C$ of realisations of $q$ we
have
\[
\dim(a_0, \dots, a_{n-1}/ A \cup C) = n.
\]
\item
For some independent sequence $(a_0, \dots, a_{n-1}, a_n)$ of realisations
of $p$ there is a countable set $C$ of realisations of $q$ such that
\[
\dim(a_0, \dots, a_{n-1}, a_n/ A \cup C) \leq n.
\]
\end{enumerate}
Then $\mathfrak{C}$ interprets a group $G$ which acts on
the geometry $P'$ induced on the realisations of $p$.
Furthermore, either $\mathfrak{C}$ interprets a non-classical group,
or $n = 1,2,3$ and
\begin{itemize}
\item
If $n = 1$, then $G$ is abelian and acts
regularly on $P'$;
\item
If $n = 2$, the action of $G$ on $P'$ is isomorphic to the affine
action of $K^+ \rtimes K^*$ on the algebraically closed field $K$.
\item
If $n = 3$, the action of $G$ on $P'$ is isomorphic to the action
of $\PGL_2(K)$ on the projective line $\mathbb{P}^1(K)$ of the
algebraically closed field $K$.
\end{itemize}
\end{theorem}

As mentioned before, 
the phrasing in terms of dimension theory is equivalent to the statement
in orthogonality calculus in Hrushovski's theorem.
The main difference with the first order result is the appearance
of the so-called {\em non-classical groups}, which 
are nonabelian $\omega$-homogeneous groups
carrying a pregeometry.
In the first order case, it follows from Reineke's theorem~\cite{Re}
that such groups cannot exist.  
Another difference is that in the interpretation,
we must use invariance rather than definability; since we have some
homogeneity in our contexts, invariant sets are definable in infinitary
logic (in the excellent case, for example, they are type-definable).

The paper is divided into four sections.
The first two sections are group-theoretic and, although motivated
by model theory, contain none. 
The first section is concerned with generalising classical theorems
on strongly minimal saturated groups and fields.
We consider groups and fields whose universe carries an $\omega$-homogeneous 
pregeometry.  We introduce generic elements and ranks, but
make no stability assumption.
We obtain a lot of information on the structure of non-classical groups, 
for example they are not solvable, 
their center is $0$-dimensional, and the quotient with the center
is divisible and torsion-free.
Nonclassical groups are very complicated; in addition to the properties
above, any two nonidentity elements of the quotient with the center 
are conjugate.
Fields carrying an $\omega$-homogeneous 
pregeometry are more amenable; as in the first order
case, we can show that they are algebraically closed.

In the second section, we generalise the theory of groups acting
on strongly minimal sets.
We consider groups $G$ {\em $n$-acting} on a pregeometry $P$,
{\em i.e.}
the action of the group $G$ respects the pregeometry, and further
(1) the integer $n$ is maximal such that
for each pair of independent $n$-tuples of the pregeometry $P$, there
exists an element of the group $G$ sending one $n$-tuple to the other, 
and (2) two elements of the group $G$ whose actions agree on an
$(n+1)$-dimensional set are identical.
As a nontriviality condition,
we require that this action must be $\omega$-homogeneous (in \cite{Hy}
Hyttinen considered this context under a stronger assumption of homogeneity,
but in order to apply the results to excellent classes we must weaken it).
We
are able to obtain a picture very similar to the
classical first order case.
We prove (see the section for precise definitions):
\begin{theorem}\label{t:intro2}
Suppose $G$ $n$-acts on a geometry $P'$. 
If $G$ admits hereditarily unique generics with respect to
the automorphism group $\Sigma$, then either there is
an $A$-invariant non-classical unbounded subgroup of $G$
(for some finite $A \subseteq P'$), or $n = 1,2,3$ and
\begin{itemize}
\item 
If $n = 1$ then $G$ is abelian and acts
regularly on $P'$.
\item 
If $n = 2$ the action of $G$ on $P'$ is
isomorphic to the affine action of $K \ltimes K^*$ on the algebraically
closed field $K$.
\item
If $n = 3$ the action of $G$ on $P'$ is isomorphic
to the action of $\PGL_2(K)$ on the projective line $\mathbb{P}^1$
of the algebraically closed field $K$.
\end{itemize}
\end{theorem}

The last two sections are completely model-theoretic.
In the third section, we consider the case of stable homogeneous
model theory, and in the fourth the excellent case.
In each case, the group we interpret is based on the automorphism
group of the monster model $\mathfrak{C}$:  
Let $p, q$ be unbounded types, say over a finite set $A$,
and assume that $p$ is quasiminimal.
Let $P = p(\mathfrak{C})$ and $Q = q(\mathfrak{C})$.
Bounded closure induces a pregeometry on $P$ and we let $P'$ be
its associated geometry.  In the stable homogeneous
case, the group we interpret is the group of permutations of
$P'$ induced by automorphisms of $\mathfrak{C}$ fixing $A \cup Q$
pointwise.
However, in the excellent case, we may not have enough homogeneity
to carry this out.
To remedy this, we consider the group $G$ of permutations of $P'$
which agree {\em locally} with automorphisms of $\mathfrak{C}$, 
{\em i.e.} a permutation $g$ of $P'$ is in $G$ if for any finite
$X \subseteq P$ and countable $C \subseteq Q$, there is an automorphism
$\sigma \in \aut(\mathfrak{C}/A \cup C)$ such that the permutation of $P'$
induced by $\sigma$ agrees with $g$ on $X$.
In each case, we show that the group $n$-acts on the geometry $P'$
in the sense of Section 2.
The interpretation in $\mathfrak{C}$ follows from the $n$-action.

Although the construction we provide for excellent classes works for
the stable homogeneous case also, for expositional reasons
we present the construction with the obvious group
in the homogeneous case first, and then present the modifications with the
less obvious group in the excellent case.

To apply Theorem~\ref{t:intro2} to $G$ and obtain 
Theorem~\ref{t:intro1}, it remains to show that
$G$ admits hereditarily unique generics with respect to some
group of automorphisms $\Sigma$.
For this, we deal with an invariant (and interpretable) subgroup of
$G$, the connected component, and deal with the group 
of automorphisms $\Sigma$ induced
by the {\em strong automorphisms}, {\em i.e.} automorphisms
preserving Lascar strong types.
Hyttinen and Shelah introduced Lascar strong types for the stable 
homogeneous case in \cite{HySh}; this is done without stability
by Buechler and Lessmann in \cite{BuLe}.
In the excellent case, this is done in detail in \cite{HyLe:2}; 
we only use the results over finite sets.

\section{Groups and fields carrying a homogeneous pregeometry}

In this section, we study algebraic structures carrying an $\omega$-homogeneous
pregeometry.  It is similar to the definition from \cite{Hy},
except that the homogeneity requirement is weaker.

\begin{definition}
An infinite model $M$ {\em carries an $\omega$-homogeneous pregeometry}
if there exists an invariant closure operator 
\[
\cl: \mathcal{P}(M) \rightarrow \mathcal{P}(M), 
\]
satisfying the axioms of a pregeometry with $\dim(M) = \|M\|$, 
and such that whenever $A \subseteq M$ is finite
and $a, b \not \in \cl(A)$,
then there is an automorphism
of $M$ preserving $\cl$, fixing $A$ pointwise, and sending $a$ to $b$. 
\end{definition}

\begin{remark}
In model-theoretic applications, the model $M$ is generally uncountable,
and $|\cl(A)| < \| M \|$, when $A$ is finite.  
Furthermore, if $a, b \not \in \cl(A)$ and $|A| < \| M \|$ one can 
often find an automorphism of $M$ fixing $\cl(A)$ pointwise, and not
just $A$.
However, we find this phrasing more natural and in non first order
contexts like excellence, $\omega_1$-homogeneity may fail.
\end{remark}

Strongly minimal $\aleph_0$-saturated groups are the simplest example
of groups carrying an $\omega$-homogeneous pregeometry.
In this case, Reineke's famous theorem~\cite{Re} 
asserts that it must be 
abelian.  Groups whose universe is a regular type are also of this form,
and when the ambient theory is stable, Poizat~\cite{Po} showed that they
are also abelian.
We are going to consider generalisations of these theorems, but
first, we need to remind the reader of some terminology.

Fix an infinite model $M$ and assume that it carries an $\omega$-homogeneous
pregeometry.  Following model-theoretic terminology,
we will say that a set $Z$ is {\em $A$-invariant}, where $A$ and $Z$ are
subsets of the model $M$, if 
any automorphism of $M$ fixing $A$ pointwise, fixes $Z$ setwise.
In particular, if $f : M^m \rightarrow M^n$ is $A$-invariant
and $\sigma$ is an automorphism of $M$ fixing $A$ pointwise,
then $f(\sigma(\bar{a})) = \sigma(f(\bar{a}))$, for any $\bar{a} \in M^m$.
We use the term {\em bounded} to mean of size less than $\|M\|$.

The $\omega$-homogeneity requirement has strong consequences.  
Obviously, any model carries the trivial pregeometry, but it is rarely
$\omega$-homogeneous; for example no group can carry a trivial 
$\omega$-homogeneous
pregeometry.
We list a few consequences of $\omega$-homogeneity which will be
used repeatedly.
First, if $Z$ is $A$-invariant, for some finite $A$, then either
$Z$ or $G \setminus Z$ is contained in $\cl(A)$ and hence 
has finite dimension (if not, choose
$x, y \not \in \cl(A)$, such that $x \in Z$ and $y \not \in Z$;
then some automorphism of $M$ fixing $A$ sends $x$ to $y$, contradicting
the invariance of $Z$).
Hence, if $Z$ is an $A$-invariant set, for some finite $A$, and has bounded
dimension,
then $Z  \subseteq  \cl(A)$.
It follows that if $a$ has bounded orbit under the automorphisms of $M$ fixing
the finite set $A$, then $a \in \cl(A)$.
This observation has the following consequence:

\begin{lemma}\label{l:fdim}
Suppose that $M$ carries an $\omega$-homogeneous pregeometry.
Let $A \subseteq M$ be finite.
Let $f: M^n \rightarrow M^m$ be an $A$-invariant function.
Then, for each $\bar{a} \in M^n$ we have
$\dim(f(\bar{a})/A) \leq \dim(\bar{a}/A)$.
\end{lemma}
\begin{proof}
Write $f = (f_0, \dots, f_m)$
with $A$-invariant $f_i : M^n \rightarrow M$,
for $i < m$.
Let $\bar{a} \in M^n$.
If $\dim(f(\bar{a}/A) > \dim(\bar{a}/A)$, then there is $i < m$
such that $f_i(\bar{a}) \not \in \cl(\bar{a} A)$.
But this is impossible since 
any automorphism $M$ fixing $A \bar{a}$
pointwise fixes $f_i(\bar{a})$. 
\end{proof}

We now introduce generic tuples.

\begin{definition}
Suppose that $M$ carries an $\omega$-homogeneous pregeometry.
A tuple $\bar{a} \in M^n$ is said to be {\em generic over $A$},
for $A \subseteq M$, if $\dim(\bar{a}/A) = n$.
\end{definition}

Since $M$ is infinite dimensional, for any finite $A \subseteq M$
and any $n < \omega$, there exists a generic $\bar{a} \in M^n$ over $A$.
Further, by $\omega$-homogeneity, if $\bar{a}, \bar{b} \in M^n$ are
both generic over the finite set $A$, then $\bar{a}$ and $\bar{b}$
are automorphic over $A$. 
This leads immediately to a proof of the following lemma.
\begin{lemma}
Suppose that $M$ carries an $\omega$-homogeneous pregeometry.
Let $A \subseteq M$ be finite and let $Z$ be
an $A$-invariant subset of $M^n$.
If $Z$ contains a generic tuple over $A$, then $Z$ contains all
generic tuples over $A$.
\end{lemma}

We now establish a few more lemmas in case when $M$ is a group $(G, \cdot)$.
Generic elements are particularly useful here.
For example, let $\bar{a} = (a_0, \dots, a_{n-1})$ and
$\bar{b}= (b_0, \dots, b_{n-1})$ belong to  $G^n$.
If $\bar{a}$ is generic over 
$A \cup \{ b_0, \dots, b_{n-1} \}$, then 
$(a_0 \dot b_0, \dots, a_{n-1} b_{n-1})$
is generic over $A$. (This follows immediately from Lemma~\ref{l:fdim}.)
When $n=1$, the next lemma asserts that if $H$ is a proper
$A$-invariant subgroup of $G$ ($A$ finite), then $H \subseteq \cl(A)$.

\begin{lemma}~\label{l:conn} 
Let $G$ be a group carrying an $\omega$-homogeneous pregeometry.
Suppose that $H$ is an $A$-invariant subgroup of $G^n$ (with
$A$ finite and $n < \omega$).
If $H$ contains a generic tuple over $A$, then $H = G^n$.
\end{lemma}
\begin{proof}
Let $(g_0, \dots, g_{n-1}) \in G^n$.
By the previous lemma, $H$ contains a generic tuple $(a_0, \dots, a_{n-1})$
over $A \cup \{ g_0, \dots, g_{n-1} \}$.
Then $(a_0 \cdot g_0, \dots, a_{n-1} \cdot g_{n-1})$ is also generic over $A$
and therefore belongs to $H$ by another application of the previous lemma.
It follows that $(g_0, \dots, g_{n-1}) \in H$.
\end{proof}

The previous lemma implies that groups carrying an 
$\omega$-homogeneous pregeometry are {\em connected}
(see the next definition).
\begin{definition}
A group $G$ is {\em connected} if it has no proper subgroup
of bounded index which is invariant over a finite set.
\end{definition}

We now introduce the {\em rank} of an invariant set.
\begin{definition}
Suppose that $M$ carries an $\omega$-homogeneous pregeometry.
Let $A \subseteq M$ be finite and let $Z$ be an $A$-invariant subset
of $M^n$.
The {\em rank of $Z$ over $A$}, written $rk(Z)$, 
is the largest $m \leq n$ such
that there is $\bar{a} \in Z$ with $\dim(\bar{a}/A) = m$. 
\end{definition}
Notice that if $Z$ is $A$-invariant and if $B$ contains $A$ is finite,
then the rank of $Z$ over $A$ is equal to the rank of $Z$ over $B$.
We will therefore omit the parameters $A$.
The next lemma is interesting also in the case where $n = 1$;
it implies that any invariant homomorphism of $G$ is either trivial
or onto.

\begin{lemma}\label{l:ker}
Let $G$ be a group carrying an $\omega$-homogeneous pregeometry.
Let $f : G^n \rightarrow G^n$ be an $A$-invariant homomorphism,
for $A \subseteq G$ finite.
Then 
\[
rk(\ker(f)) + rk(\ran(f)) = n.
\]
\end{lemma}
\begin{proof}
Let $k \leq n$ such that $rk(\ker(f)) = k$.
Fix $\bar{a} = (a_0, \dots, a_{n-1}) \in \ker(f)$ be of dimension
$k$ over $A$.
By a permutation, we may assume that $(a_0, \dots, a_{k-1})$
is independent over $A$.

Notice that by $\omega$-homogeneity and $A$-invariance
of $\ker(f)$, for each generic $(a'_0, \dots, a'_{i-1})$
over $A$ (for $i < k$), 
there exists $(b_0, \dots, b_{n-1}) \in \ker(f)$ such that 
$b_i = a'_i$ for $i < k$.
We now claim that for any $i < k$ and any $b \not \in \cl(A)$,
there is $\bar{b} = (b_0, \dots, b_{n-1}) \in \ker(f)$ such that $b_j = 1$
for $j < i$ and $b_i = b$.
To see this, notice that $(a_0^{-1}, \dots, a_{i-1}^{-1})$
is generic over $A$
(by Lemma~\ref{l:fdim}).  Choose $c \in G$ generic over $A \bar{a}$.
Then there is $(d_0, \dots, d_{n-1}) \in \ker(f)$ such
that $d_j = a_j^{-1}$ for $j < i$ and $d_i = c$.
Let $(e_0, \dots, e_{n-1}) \in \ker(f)$ be the product of $\bar{a}$ with
$(d_0, \dots, d_{n-1})$.
Then $e_j = 1$ if $j < i$ and $e_i = a_i \cdot c^{-1} \not \in \cl(A)$.
By $\omega$-homogeneity, there is an automorphism of $G$
fixing $A$ sending $e_i$ to $b$.
The image of $(e_0, \dots, e_{n-1})$ under this automorphism
is the desired
$\bar{b}$. 

We now show that $rk(\ran(f)) \leq n - k$.
Let $\bar{d} = f(\bar{c})$.
Observe that by multiplying $\bar{c} = (c_0, \dots, c_{n-1})$ 
by appropriate elements in $\ker(f)$,
we may assume that $c_i \in \cl(A)$ for each $i < k$.
Hence $\dim(\bar{c}/A) \leq n-k$ so the conclusion follows
from Lemma~\ref{l:fdim}.

To see that $rk(\ran(f)) \geq n-k$, choose $\bar{c} \in G^n$
such that $c_i = 1$ for $i < k$ and $(c_k, \dots , c_{n-1})$
is generic over $A$.
It is enough to show that $\dim(f(\bar{c})/A) \geq \dim(\bar{c}/A)$.
Suppose, for a contradiction, that 
$\dim(f(\bar{c})/A) < \dim(\bar{c}/A)$.
Then there is $i < n$, with $k \leq i$ such that
$c_i \not \in \cl(f(\bar{c}) A)$.
Let $d \in G \setminus \cl(A f(\bar{c}) \bar{c})$ and
choose an automorphism $\sigma$ fixing $A f(\bar{c})$ such
that $\sigma(c_i) = d$.
Let $\bar{d} = \sigma(\bar{c})$. 
Then
\[
f(\bar{d}) = f(\sigma(\bar{c})) = \sigma( f(\bar{c})) = f(\bar{c}).
\]
Let $\bar{e} = (e_0, \dots, e_{n-1}) = \bar{c} \cdot \bar{d} ^{-1}$.
Then $\bar{e} \in \ker(f)$, $e_j = 1$ for $j < k$, and 
$e_i = c_i \cdot d^{-1} \not \in \cl(A)$.
By $\omega$-homogeneity, we may assume that $e_i
\not \in \cl(A \bar{a})$.
But $\bar{a} \cdot \bar{e} \in \ker(f)$,
and $\dim(\bar{a} \cdot \bar{e}/A) \geq k+1$ 
(since the $i$-th coordinate of $\bar{a} \cdot \bar{e}$ 
is not in $\cl(a_0, \dots, a_{k-1} A)$).
This contradicts the assumption that $rk(\ker(f)) = k$.
\end{proof}

The next theorem is obtained by adapting Reineke's proof to our
context. 
For expository purposes, 
we sketch some of the proof and refer to reader to \cite{Hy} for details.
We are unable to conclude that groups carrying an 
$\omega$-homogeneous pregeometry
are abelian, but we can still obtain a lot of information.

\begin{theorem}\label{t:ugly} Let $G$ be a nonabelian group which carries
an $\omega$-homogeneous pregeometry. 
Then the center $Z(G)$ has dimension $0$, $G$ is not solvable, 
any two nonidentity elements in the quotient group $G/Z(G)$ 
are conjugate, and $G/Z(G)$ is torsion-free and divisible.
Also $G$ contains a free subgroup on $\dim(G)$ many generators, 
and the first order theory of $G$ is unstable.
\end{theorem}

\begin{proof}
If $G$ is not abelian, then the center of $G$,
written $Z(G)$, is a proper subgroup of $G$.
Since $Z(G)$ is invariant, Lemma~\ref{l:conn} implies that
$Z(G) \subseteq \cl(\emptyset)$.
 
We now claim that if $H$ is an $A$-invariant proper normal subgroup of
$G$ then $H \subseteq Z(G)$.
 
By the previous lemma, $H$ is finite dimensional.
For $g, h \in H$, define 
$X_{g,h} = \{ x \in  G : g^x = h \}$.
Suppose, for a contradiction, that $H \not \subseteq Z(G)$ and
choose $h_0 \in H \setminus Z(G)$.
If for each $h \in H$, the set $X_{h_0,h}$ is finite dimensional,
then $X_{h_0,h_1} \subseteq \cl(h_0 h)$, and so
$G \subseteq \bigcup_{h \in H} X_{h_0,h} \subseteq \cl(H)$,
which is impossible since $H$ has finite dimension.
Hence, there is $h_1 \in H$, 
such $X_{h_0,h_1}$ is infinite dimensional and has
finite-dimensional complement.
Similarly, there is $h_2 \in H$ such that $X_{h_1, h_2}$
has finite dimensional complement.
This allows us to choose $a, b \in G$ such that
$a, b, ab$ belong to both $X_{h_0,h_1}$ and $X_{h_1,h_2}$.
Then, $h_1 = h_0^{ab}=(h_0^b)^a = h_2$.
This implies that the centraliser of $h_1$ has infinite dimension
(since it is $X_{h_1,h_2}$) and must therefore be all of $G$ 
by the first paragraph of this proof.
Thus $h_1 \in Z(G)$, which is impossible, since it is the conjugate
of $h_0$ which is not in $Z(G)$.

We now claim that $G^* = G/Z(G)$ is not abelian.
Suppose, for a contradiction, that $G^*$ is abelian. 
Let $a \in G \setminus Z(G)$.
Then the sets $X_{a,k} = \{ b \in G:  a^b = ak \}$, where
$k \in Z(G)$ form a partition of $G$, and so, as above, 
there is $k_a \in Z(G)$ such that 
$X_{a, k_a}$ has finite-dimensional complement.
Now notice that $k_a = 1$: Otherwise, for $c \in X_{a, k_a}$, 
the infinite dimensional sets $cX_{a, k_a}$ and $X_{a, k_a}$ are disjoint,
which is impossible, since they each have finite-dimensional complements.
But then, $X_{a, 1}$ is a subgroup of $G$ of infinite
dimension
and so is equal to $G$, which implies that $a \in Z(G)$, a contradiction.

Since $G / [G,G]$ is abelian, and $[G,G]$ is normal and invariant,
then it cannot be proper (otherwise $[G,G] \leq Z(G)$).
It follows that $G$ is not solvable.

It follows easily from the previous claims that $G^*$ is centerless.  
We now show that any two nonidentity elements in $G^*$ are conjugate:
Let $a^* \in G^*$ be a nonidentity element.
Since $G^*$ is centerless, the centraliser of $a^*$ in $G^*$ is
a proper subgroup of $G^*$.  Hence, the inverse image of this centraliser
under the canonical homomorphism induces a proper subgroup of $G$,
which must therefore be of finite dimension. 
Hence, the set of conjugates of $a^*$ in $G^*$ is all of $G^*$,
except for a set of finite dimension.
It then follows that the set of elements of $G^*$ which are not 
conjugates of $a^*$ must have bounded dimension.

Since this holds for any nonidentity $b^* \in G^*$, this implies
that any two nonidentity elements of $G^*$ must be conjugates.
The instability of $Th(G)$ now follows as in the proof of \cite{Po}:
Since $G/Z(G)$ is not abelian and any two nonidentity elements of it are
conjugate, we can construct an infinite strictly ascending chain
of centralisers.  This contradicts first order stability.

That $G$ is torsion free and divisible  
is proved similarly (see \cite{Hy}
for details). Finally, it is easy to check that any independent subset
of $G$ must generate a free group.
\end{proof}

Hyttinen called such groups {\em bad} in \cite{Hy}, but this conflicts with
a standard notion, so we re-baptise them:
\begin{definition}
We say that a group $G$ is {\em non-classical} if it is nonabelian
and carries
an $\omega$-homogeneous pregeometry.
\end{definition}

\begin{question}
Are there non-classical groups?  And if there are, can they arise in the 
model-theoretic contexts we consider in this paper?
\end{question}

We now turn to fields.
Here, we are able to adapt the proof of Macintyre's classical
theorem~\cite{Ma} that $\omega$-stable fields are algebraically closed.

\begin{theorem}\label{t:ma}
A field carrying an $\omega$-homogeneous pregeometry
is algebraically closed.
\end{theorem}
\begin{proof}
To show that $F$ is algebraically closed, it is enough to show
that any finite dimensional field extension $K$ of $F$ is perfect,
and has no Artin-Schreier or Kummer extension.

Let $K$ be a field extension of $F$  
of finite degree $m < \omega$.
Let $P \in F[X]$ be an irreducible polynomial of degree $m$
such that $K = F(\xi)$, where $P(\xi) = 0$.
Let $A$ be the finite subset of $F$ consisting of the coefficients
of $P$. 
We can represent $K$ in $F$ as follows:  $K^+$ is the
vector space $F^m$, {\em i.e.}
$a = a_0 + a_1 \xi + \dots a_{m-1} \xi^{m-1}$ is represented
as $(a_0, \dots, a_{m-1})$.
We can then easily represent addition in $K$ and
multiplication (the field product in $K$ induces a bilinear form on $(F^+)^m$)
as $A$-invariant operations. 
Notice that an automorphism $\sigma$ of $F$ fixing $A$ pointwise
induces an automorphism of $K$, via
\[
(a_0, \dots , a_{m-1}) \mapsto (\sigma(a_0), \dots, \sigma(a_{m-1})).
\]

We now consider generic elements of the field.
For a finite subset $X \subseteq F$ containing $A$,
we say that an $a \in K$
is {\em generic over $X$} if $\dim(a_0 \dots a_{m-1}/X) = m$
(that is $(a_0, \dots, a_{m-1})$ is generic over $X$),
where 
$a_i \in F$ and 
$a = a_0 + a_1 \xi + \dots a_{m-1} \xi^{m-1}$.
Notice that if $a, b \in K$ are generic over $X$ (with $X \subseteq
F$ finite containing $A$) then there exists an automorphism of $K$
fixing $X$ sending $a$ to $b$.
We prove two claims about generic elements.
\begin{claim}
Assume that $a \in K$ is generic over the finite set $X$, 
with $A \subseteq X \subseteq F$.
Then $a^n$, $a^n-a$ for $n < \omega$, as well as $a+b$ and $ab$ for 
$b= b_0 + b_1 \xi + \dots b_{m-1} \xi^{m-1}$, $b_i \in X$ ($i < m$) are
also generic over $X$.
\end{claim}  
\begin{proof}[Proof of the claim]
We prove that $a^n$ is generic over $X$. The other proofs are similar.
Suppose, for a contradiction,
that $a^n = c_0 + c_1 \xi +  \dots + c_{m-1} \xi^{m-1}$
and $\dim(c_0 \dots c_{m-1}/X) < m$.
Then $\dim(a_0 \dots a_{m-1}/Xc_0 \dots c_{m-1}) \geq 1$,
so there is $a_i \not \in \cl(Xc_0 \dots c_{m-1})$.
Since $F$ is infinite dimensional, there are infinitely many
$b \in F \setminus \cl(Xc_0 \dots c_{m-1})$, and
by $\omega$-homogeneity there is an automorphism of $F$
fixing $X c_0 \dots c_{m-1}$ sending $a_i$ to $b$.
It follows that there are infinitely many $x \in K$ such
that $x^n = a^n$, a contradiction.
\end{proof}

\begin{claim}
Let $G$ be an $A$-invariant subgroup of $K^+$ (resp. of $K^*$).
If $G$ contains an element generic over $A$ then $G = K^+$ (resp. $G = K^*$).
\end{claim}
\begin{proof}[Proof of the claim]
We prove only one of the claims, as the other is similar.
First, observe that if $G$ contains an element of $K$ generic over $A$,
it contains all elements of $K$ generic over $A$.
Let $a \in K$ be arbitrary.
Choose $b \in K$ generic over $Aa$.  Then $b \in G$, and 
since 
$a + b$ is generic over $Aa$ (and hence over $A$), we have 
also $a+b \in G$.  It follows that $a \in G$, since $G$ is a
subgroup of $K^+$. Hence $G = K^+$. 
\end{proof}

Consider the $A$-invariant subgroup $\{ a^n : a \in K^* \}$ of $K^*$.
Let $a \in K$ be generic over $A$.  Since $a^n$ is generic over $A$
by the first claim, we have that $\{ a^n : a \in K^* \} = K^*$
by the second claim.
This shows that $K$ is perfect (if the characteristics is a prime $p$,
this follows from the existence of $p$-th roots, and every field
of characteristics $0$ is perfect).

Suppose $F$ has characteristics $p$.
The $A$-invariant subgroup $\{ a^p - a : a \in K^+ \}$ of $K^+$ contains
a generic element over $A$ and hence 
$\{ a^p - a : a \in K^+ \} = K^+$.

The two previous paragraphs show that $K$ is perfect and has no Kummer
extensions (these are obtained by adjoining a solution to
the equation $x^n = a$, for some $a \in K$)
or Artin-Schreier extensions (these are obtained by adjoining a solution 
to the equation
$x^p - x = a$, for some $a \in K$, where $p$ is the characteristics).
This finishes the proof.
\end{proof}

\begin{question}
If there are non-classical groups, are there also division rings
carrying an $\omega$-homogeneous pregeometry which are not fields?
\end{question}

\section{Group acting on pregeometries}

In this section, we generalise some classical results
on groups acting on strongly minimal sets. 
We recall some of the facts, terminology, and results  
from \cite{Hy}, and then prove some additional theorems.

The main concept is that of {\em $\Sigma$-homogeneous group $n$-action
of a group $G$ on a pregeometry $P$}. This consists of the following:

We have a group $G$ acting on a pregeometry $(P, \cl)$.
We denote by $\dim$ the dimension inside the pregeometry $(P, \cl)$
and always assume that $\cl(\emptyset) = \emptyset$.
We write 
\[
(g, x) \mapsto gx, 
\]
(or sometimes $g(x)$ for legibility) 
for the action of $G$ on $P$.
For a tuple $\bar{x} = (x_i)_{i<n}$ of elements of $P$, we 
write $g\bar{x}$ or $g(\bar{x})$ for $(gx_i)_{i<n}$.
The group $G$ acts on the universe of $P$ and respects the pregeometry,
{\em i.e.} 
\[
a \in \cl(A)
\quad 
\text{if and only if}
\quad
ga \in \cl(g(A)),
\]
for $a \in P$,
$A \subseteq P$ and $g \in G$.

We assume that the action of $G$ on $P$ is an {\em $n$-action},
{\em i.e.} has the following two properties:
\begin{itemize}
\item
The action has {\em rank $n$}:
Whenever $\bar{x}$ and $\bar{y}$ are two $n$-tuples of elements of
$P$ such that $\dim(\bar{x}\bar{y}) = 2n$, 
then there is $g \in G$ such that $g\bar{x}=\bar{y}$.
However, for some $(n+1)$-tuples $\bar{x}, \bar{y}$ 
with $\dim(\bar{x} \bar{y}) = 2n+2$,
there is no $g \in G$ is such that $g\bar{x}=\bar{y}$.

\item
The action is {\em $(n+1)$-determined}:
Whenever the action of $g, h \in G$ agree on an $(n+1)$-dimensional
subset $X$ of $P$, then $g=h$.
\end{itemize}

An {\em automorphism of the group action}
is a pair of automorphisms $(\sigma_1, \sigma_2)$,
where $\sigma_1$ is an automorphism of the group $G$ and
$\sigma_2$ is an automorphism of the pregeometry $(P, \cl)$,
which preserve the group action, 
{\em i.e.} 
\[
\sigma_2(gx) = \sigma_1(g) \sigma_2(x).
\]
Following model-theoretic practice, we will simply think of 
$(\sigma_1, \sigma_2)$ as a single automorphism $\sigma$
acting on two disjoint structures (the group and the pregeometry)
and write $\sigma(gx) = \sigma(g) \sigma(x)$.

We let $\Sigma$ be a group of automorphisms of this group action.
We assume that the group action is 
{\em $\omega$-homogeneous with respect to $\Sigma$}, {\em i.e.}
if whenever $X \subseteq P$ is 
{\em finite} and $x, y \in P \setminus \cl(X)$, then there is 
an automorphism $\sigma \in \Sigma$ such
that $\sigma(x) = y$ and $\sigma \restriction X = id_X$.
Notice that $x \in P$ is fixed under all automorphisms in $\Sigma$
fixing the finite set $X$ pointwise, then $x \in \cl(X)$.

This is essentially the notion that Hyttinen isolated in
\cite{Hy}.  There are two slight differences:  
(1) We specify the automorphism
group $\Sigma$, whereas \cite{Hy} works with {\em all} automorphisms
of the action (but there he allows extra structure on $P$,
thus changing the automorphism group, 
so the settings are equivalent).
(2) We require the existence of
$\sigma \in \Sigma$ such
that $\sigma(x) = y$ and $\sigma \restriction X = id_X$,
when $x, y \not \in \cl(X)$ only for {\em finite} $X$.
All the statements and proofs from \cite{Hy} can be easily modified.
Some of the results of this section are easy adaptation from the proofs
in \cite{Hy}.  To avoid unnecessary repetitions, we sometimes list
some of these results as facts and refer the reader to \cite{Hy}.

Homogeneity is a nontriviality condition;
it actually has strong consequences. For example, 
if $\bar{x}$ and $\bar{y}$ are $n$-tuples each of dimension $n$, then
there is $g \in G$ such that $g(\bar{x}) = \bar{y}$.
Further, for no pair
of $(n+1)$-tuples $\bar{x}$, $\bar{y}$ with $\dim(\bar{x}\bar{y}) = 2n+2$ 
is there a $g \in G$ sending $\bar{x}$ to $\bar{y}$.
This implies that 
if $\bar{x}$ is an independent $n$-tuple and $y$ is an element outside
$\cl(\bar{x}g(\bar{x}))$,
then necessarily $g(y) \in \cl(\bar{x}g(\bar{x}) y)$.

We often just talk about {\em $\omega$-homogeneous group
acting on a pregeometry}, when the identity
of $\Sigma$ or $n$ are clear from the context.

The classical example of homogeneous group actions on a pregeometry
are  
definable groups acting on a strongly minimal sets inside a saturated
model.  Model theory provides important tools to deal with this situation;
we now give generalisations of these tools and define 
{\em types}, {\em stationarity}, {\em generic elements}, {\em connected
component}, and so forth in this general context.

\relax From now until the end of this section, we fix
an $n$-action of $G$ on
the pregeometry $P$ which is $\omega$-homogeneous with respect to $\Sigma$.

Let $A$ be a $k$-subset of $P$ with $k < n$.
We can form a new homogeneous group action by {\em localising at $A$}:
The group $G_A \subseteq G$ is the stabiliser of $A$; 
the pregeometry
$P_A$ is obtained from $P$ by considering the new closure
operator $\cl_A(X) = \cl(A \cup X) \setminus \cl(A)$
on $P \setminus \cl(A)$; 
the action of $G_A$ on $P_A$ 
is by restriction;
and let $\Sigma_A$ be the group
of automorphisms in $\Sigma$ fixing $A$ pointwise.
We then have a $\Sigma_A$-homogeneous group $(n-k)$-action of $G_A$ on
the pregeometry $P_A$.

Generally, for $A \subseteq G \cup P$, we denote by 
$\Sigma_A$ the group of automorphisms in $\Sigma$ which fix
$A$ pointwise.

Using $\Sigma$, we can talk about {\em types} of elements of
$G$: these are the orbits of elements of $G$ under $\Sigma$.
Similarly, the {\em type of an element $g \in G$ over $X \subseteq P$}
is the orbit of $g$ under $\Sigma_X$.
We write $\tp(g/X)$
for the type of $g$ over $X$.

\begin{definition}
We say that $g \in G$ is {\em generic} over $X \subseteq P$,
if there exists an independent $n$-tuple $\bar{x}$ of $P$ 
such that
\[
\dim(\bar{x} g(\bar{x}) / X) = 2n.
\]
\end{definition}

It is immediate that if $g$ is generic over $X$ then so is its inverse.
An important property is that
given a finite set $X \subseteq P$,
there is a $g \in G$ generic over $X$.
Notice also that genericity of $g$ over $X$ is a property
of $\tp(g/X)$; we can therefore talk about {\em generic types over $X$},
which are simply types of elements generic over $X$.
Finally, if $\tp(g/X)$ is generic over $X$, $X \subseteq Y$ are finite
dimensional, then there is $h \in G$ generic over $Y$
such that $\tp(h/X) = \tp(g/X)$.

We can now define stationarity in the natural way (notice the extra
condition on the number of types; this condition holds trivially
in model-theoretic contexts).
\begin{definition}
We say that $G$ is {\em stationary} with respect to $\Sigma$, 
if whenever $g, h \in G$ with $\tp(g/\emptyset) = \tp(h/\emptyset)$
and $X \subseteq P$ is finite and both $g$ and $h$ are generic
over $X$, then $\tp(g/X) = \tp(h/X)$.  
Furthermore, we assume that 
the number of types over each finite set is bounded.
\end{definition}

The following is a strengthening of stationarity.
\begin{definition}
We say that $G$ has {\em unique generics} if for all finite
$X \subseteq P$ and $g, h \in G$ generic over $X$ we have
$\tp(g/X) = \tp(h/X)$.
\end{definition}

We now introduce the {\em connected component $G^0$}:
We let $G^0$ be the intersection of all invariant, normal subgroups
of $G$ with bounded index. Recall that a set is {\em invariant}
(or more generally {\em $A$-invariant}) if it is fixed setwise
by any automorphism in $\Sigma$ ($\Sigma_A$ respectively).
The proof of the next fact is left to the reader; it can also
be found in \cite{Hy}.
\begin{fact}\label{f:connected}
If $G$ is stationary, then $G^0$ is a normal invariant
subgroup of $G$ of bounded index.
The restriction of the action of $G$ on $P$ to $G^0$ is an $n$-action,
which is homogeneous with respect to the group
of automorphisms obtained from $\Sigma$ by restriction.
\end{fact}

We provide the proof of the next proposition to convey the flavour of these
arguments.

\begin{proposition}\label{t:unique}
If $G$ is stationary
then $G^0$ has unique generics.
\end{proposition}

\begin{proof}
Let $Q$ be the set of generic types over the empty set.
For $q \in Q$ and $g \in G$, we define $gq$ as follows:
Let $X \subseteq P$ with the property that
$\sigma \restriction X = id_X$ implies $\sigma(g) = g$ for
any $\sigma \in \Sigma$.
Choose $h \models q$ which
is generic over $X$.  Define $gq = \tp(gh/\emptyset)$.

Notice that by stationarity of $G$, the definition of $gq$
does not
depend on the choice of $X$ or the choice of $h$.
Similarly, the value of $gq$ depends on $\tp(g/\emptyset)$ only.
We claim that 
\[
q \mapsto gq
\] 
is a group action of $G$ on $Q$.
Since $1q = q$, in order to prove that this is indeed an action on $Q$, 
we need to show that $gq$ is generic and $(gh)(q) = g(hq)$.

This is implied by the following claim:
If $X  \subseteq P$ is finite 
containing $\bar{x}$ and $g(\bar{x})$, 
where $\bar{x}$ is an independent $(n+1)$-tuple
of elements in $P$,
and $h \models q$ is generic over $X$, then $gh$ is generic over $X$.

To see the claim, choose $\bar{z}$ an $n$-tuple of elements of $P$ such that
\[
\dim(\bar{z} h(\bar{z}) /X) = 2n.
\]
Notice that
$h(\bar{z}) \subseteq \cl(X gh(\bar{z}))$, 
since any $\sigma \in \Sigma$ fixing $X gh(\bar{z})$
pointwise fixes $h(\bar{z})$ (for any such $\sigma$, we have
$\sigma(h(\bar{z})) = \sigma (g^{-1} gh(\bar{z}))=
\sigma (g^{-1})\sigma (gh(\bar{z})) = g^{-1}gh(\bar{z}) = h(\bar{z})$).
Thus, $\dim(\bar{z} gh(\bar{z}) /X) \geq \dim(\bar{z}h(\bar{z})/X) = 2n$, 
so $\bar{z}$
demonstrates that $gh$ is generic over $X$.

Now consider the kernel $H$ of the action, namely the
set of $h \in G$ such that $hq = q$ for each $q \in Q$.
This is clearly an invariant subgroup, and since the action depends
only on $\tp(h/\emptyset)$, $H$ must have bounded index
(this condition is part of the definition of stationarity).
Hence, by definition, the connected component $G^0$ is a subgroup of $H$.

By stationarity of $G$, if $G^0$ does not have unique generics, there
are $g, h \in G^0$ be generic over the empty set
such that $\tp(g/\empty) \not = \tp(h/\empty)$.
Without loss of generality, we may assume that $h$ is generic over 
$\bar{x} g(\bar{x})$, 
where $\bar{x}$ is an independent
$(n+1)$-tuple of $P$.
Now it is easy to check that
$hg^{-1}(\tp(g/\emptyset)) = \tp(h /\emptyset))$, so that
$hg^{-1} \not \in H$.  But $hg^{-1}h \in G^0 \subseteq H$, 
since $g,h \in  G^0$,
a contradiction.
\end{proof}

We now make another definition:
\begin{definition}
We say that 
$G$ {\em admits hereditarily unique generics} if
$G$ has unique generics and
for any independent $k$-set $A \subseteq P$ with $k < n$,
there is a normal subgroup $G'$ of $G_A$ such that
the action of $G'$ on $P_A$ is a homogeneous $(n-k)$-action 
which has {\em unique generics}. 
\end{definition}

If we have a $\Sigma$-homogeneous $1$-action of a group $G$
on a pregeometry $P$ which has unique generics,
then 
the pregeometry lifts up on the universe of the group in the natural
way and so the group carries a homogeneous pregeometry:
For $g \in G$ and $g_0, \dots, g_k \in G$,
we let
\[
g \in \cl(g_0, \dots, g_k),
\]
if for some independent $2$-tuple $\bar{y} \in P$ and some $x \in P \setminus
\cl(\bar{y} g(\bar{y}) g_0(\bar{y}) \dots g_k(\bar{y}))$ then
\[
g(x) \in \cl(x g_0(x), \dots, g_k(x)).
\]

Notice first that this definition does not depend on the choice
of $x$ and $\bar{y}$: 
Let $x' \not \in \cl(\bar{y}' g(\bar{y}') g_0(\bar{y}') \dots g_k(\bar{y}'))$
for another independent $2$-tuple $\bar{y}'$.
Let $z$ be such that 
\[
z \not \in \cl(\bar{y} g(\bar{y}) g_0(\bar{y}) \dots 
g_k(\bar{y})\bar{y}' g_0(\bar{y}') \dots g_k(\bar{y}')).
\]
Then by homogeneity, there exists 
$\sigma \in \Sigma_{\bar{y} g_0(\bar{y}) \dots g_k(\bar{y})}$
such that $\sigma(x) = z$, and 
$\tau \in \Sigma_{\bar{y}' g_0(\bar{y}') \dots g_k(\bar{y}')}$
such that $\tau(z) = x'$.  Notice that 
$\sigma(g) = \tau(g) = g$ and $\sigma(g_i) = \tau(g_i) = g_i$
for $i \leq k$ by $2$-determinacy.
Hence $g(x) \in \cl(x g_0(x), \dots ,g_k(x))$ if and only
if $g(x') \in \cl(x' g_0(x), \dots ,g_k(x))$ by applying $\sigma \circ \tau$.

We define $g \in \cl(A)$ for $g, A$ in $G$, where $A$ may be infinite, 
if there are $g_0, \dots, g_k\in G$
such that $g \in \cl(g_0, \dots, g_k)$.  It is not difficult to check
that this
induces a pregeometry on $G$.

The unicity of generics implies that the pregeometry is 
$\omega$-homogeneous:
Suppose $g, h \not \in \cl(A)$,
where $A \subseteq G$ is finite.
For a tuple $\bar{z}$, write $A(\bar{z}) = \{ f(\bar{z}) : f \in A\}$.
Let $\bar{y}$ be an independent $2$-tuple
and choose $x \not \in \cl(\bar{y}g(\bar{y}) h(\bar{y}) A(\bar{y}))$
with $g(x), h(x) \not \in \cl(x A(x))$.
Since $G$ has unique generics, it is enough to show that
$g,h$ are generic over $\bar{y}A(\bar{y})$.
Let $z\in P$ outside $\cl(xg(x) h(x) A(x))$.
Then, since the action has rank $1$, we must have
$f(x) \in \cl(x A(x))$, for each $f \in A$.
Hence $\cl(xz A(x) A(z) \subseteq \cl(x z A(x))$
and by exchange, this implies that 
$g(x), h(x) \not \in \cl(x z A(x))$.
Let $z'$ be an element outside $\cl(x z A(x) g(x) h(x))$.
It is easy to see that $\dim(z' g(z') / x z A(x) A(z)) = 2$
and so $g$ is generic over $x z A(x) A(z)$ and hence over $x A(x)$.
The same argument shows that $h$ is generic over $x A(x)$.
Hence, there is $\sigma \in \Sigma$ fixing $A$ such that $\sigma(g) = h$.

We have just proved the following fact:
\begin{fact}
If $n = 1$, $G$ is stationary and has unique generics, 
then $G$ carries an $\omega$-homogeneous
pregeometry.
\end{fact}

Admitting hereditarily unique generics is connected to $n$-determinacy
and non-classical groups in the following way.
The proof of the next fact is in \cite{Hy}; notice that the group
$(G_A)^0$ $1$-acts and so carries an $\omega$-homogeneous pregeometry.

\begin{fact}\label{f:ndet}
Suppose that $G$ admits hereditarily unique generics.
Then either $(G_A)^0$ is non-classical, 
for some independent $(n-1)$-subset $A \subseteq P$ 
or the action of $G$ on $P$ is $n$-determined.
\end{fact}

So in the case of $n = 1$, either the connected component is non-classical, 
or it is abelian and the action of $G$ on $P$ is $1$-determined.
Hence, the action of $G^0$ on $P$ is regular.

Again, see \cite{Hy} for the next fact.
\begin{fact}\label{f:3}
If the action is $n$-determined 
then $n = 1, 2, 3$.
\end{fact}

Following standard terminology, we set:
\begin{definition}
We say that the $n$-action of $G$ on $P$ is {\em sharp}
if it is $n$-determined.
\end{definition}
Notice that if $G$ $n$-acts sharply on $P$, then
the element of $G$ sending a given independent $n$-tuple
of $P$ to another is unique.

\relax From now, until theorem~\ref{t:field} we assume that the
$n$-action of $G$ on $P$ is
sharp.
Hence $n = 1, 2, 3$ by Fact~\ref{f:3}.
We are interested in constructing a field so we may assume that $n \geq 2$.
By considering the group $G_a$ acting on the pregeometry $P_a$ with
$\Sigma_a$ when $a$ is an element of $P$, we may assume that
$n =2$.
This part has not been done in \cite{Hy}.

Following Hrushovski~\cite{Hr}, 
we now introduce some invariant subsets of $G$, which will
be useful in the construction of the field.
We first consider the set of involutions.
\begin{definition}
Let $I = \{ g \in G : g^2 = 1 \}$.
\end{definition}

The set $I$ may not be a group.
\begin{definition}
Let $a \in P$.
We let $N_a \subseteq G$ consists of those elements $g \in G$
for which the set  
\[
\{ h(a) : h \in I, gh \not \in I \}
\]
has bounded dimension in $P$.
\end{definition}

We now establish a few facts about $I$ and $N_a$; in particular
that $N_a$ is an abelian subgroup of $G$:
\begin{lemma}~\label{l:31} Let $a \in P$.
\begin{enumerate}
\item
Let $g, h \in I$.  If $g(a) = h(a)$ and $g(a) \not \in \cl(a)$, then 
$g =h$.
\item
Let $g, h \in I$.
Assume that $g(a) \not \in \cl(a h(a))$, and $h(a) \not \in \cl(a
g(a))$.
Then $g h \in N_a$.
\item
Let $g, h \in N_a$.
If $g(a) = h(a)$, then $g = h$.
\item
$N_a$ is a subgroup of $G$.
\end{enumerate}
\end{lemma}

\begin{proof} (1)  Since $g^2 = h^2 = 1$, then $g(g(a)) = a$,
and $h(g(a)) = a$, since $g(a) = h(a)$.
Hence $g$ and $h$ agree on a 2-dimensional set so $g = h$ 
since the action of $G$ is 2-determined.

(2)
It is easy to see that $h(a) \not \in \cl(a gh(a))$.
Now $ghh=g \in I$, since both $g, h \in I$.  
But then, $ghf \in I$ for all generic $f \in I$.
Hence, $gh \in N_a$.

(3) Suppose first that $a \not \in \cl(g(a))$. 
Choose $f \in I$ and $b \in P$ such that $b \not \in \cl(a g(a))$
and $f(b) = a$.
Then $gf$ and $hf$ belong to $I$ and since $gf(b) = hf(b)$,
we have $gf = hf$ by (1) so $g = h$. 

Now if $g(a) = a$, we show that $g = 1$.
If not, then since the action is $2$-determined we have that $g(b) \not = b$,
for any $b \in P$ with $b \not \in \cl(a)$.
Now let $f \in I$ be such that $f(a) = b$ for $b \not \in \cl(a)$.
Then $gfgf(a) = a$, since $g \in N_a$. 
But this implies that $g(b) = b$, a contradiction.

(4)  Let $g, h \in N_a$:  First we show that $gh \in N_a$.
Choose $f \in I$ such that $f(a) \not \in \cl(a g(a) h(a))$.
Then $h(f(a)) \not \in \cl(a g(a))$.
Hence, since $h \in N_a$ we have that $hf \in I$, so $ghf \in I$
since $g \in N_a$. This shows that $gh \in N$.
Second, we show that $g^{-1} \in N_a$.  If $g^2 = 1$, then it is clear.
Otherwise by (3) $g(a) \not \in \cl(a)$.
Let $f \in I$ such that $f(a) \not \in \cl(a g(a))$.
Then $gf \in I$ and so $gfgf = 1$ so that $g^{-1} = fgf$.
But, by (2) $fgf \in N_a$.
\end{proof}

\begin{lemma} 
For $a \in P$ the group $N_a$ is abelian.
\end{lemma}
\begin{proof}
By $2$-determinacy and Lemma~\ref{l:31}, it is easy to
verify that $N_a$ carries a homogeneous pregeometry
$(N_a, \cl')$:  For $X \subseteq N_a$,
let $g \in \cl'(X)$ if $g(a) \in \cl(a \cup \{f(a) : f \in X \})$.
It is $\omega$-homogeneous with respect to the restrictions
of $\sigma \in \Sigma_a$ to $N_a$.
If $N_a$ were not abelian, then its center $Z(N_a)$ be 
$0$-dimensional 
and using Lemma~\ref{l:31}~(3) it follows that $Z(N_a)$ is trivial.
Also there is $g \in N_a$ with $g \not = g^{-1}$. 
By Theorem~\ref{t:ugly}, choose $f \in G$ such that $g^{-1} = f^{-1}gf$.
Let $h \in I$ be independent from $g$ and $f$ (in the
sense of the pregeometry $\cl'$ on $N_a$), and as in the proof
of Lemma~\ref{l:31} we have $g = hg^{-1} h^{-1}= hf^{-1}gfh$.
Then $fh \in I$ and since $hf^{-1}=(fh)^{-1}=fh$, there
is $k \in I$ independent from $g$ such that $g = kgk$.
But $kgk= g^{-1}$, a contradiction.
\end{proof}

We can now state a proposition.
Recall that we say that a group action is {\em regular}
if it is sharply transitive.

Recall that a {\em geometry} is a pregeometry such that 
$\cl(a) = \{ a \}$ for each $a \in P$ (we already assumed
that $\cl(\emptyset) = \emptyset$.

\begin{proposition}~\label{t:field}
Consider the $\Sigma$-homogeneous sharp $2$-action of $G$ on $P$.
Then, $G_a$ acts regularly on $N_a$ by conjugation and 
$G = G_a \ltimes N_a$.
Furthermore, either $G_a$ is non-classical, or $G_a$ is abelian
and the action of $G_a$ on $N_a$
induces the structure of an algebraically closed field on $N_a$. 
Furthermore, if $G_a$ is abelian, then the action of $G$ on $P$
is sharply $2$-transitive (on the {\em set} $P$), and $P$ is a geometry.  
\end{proposition}
\begin{proof}
We have already shown that $N_a$ carries a homogeneous pregeometry, and
$G_a$ carries a homogeneous pregeometry by $1$-action.
Furthermore, $N_a$ is abelian.

We now show that $G_a$ acts on $N_a \setminus \{0 \}$ by conjugation, {\em
i.e.} if $g \in N_a$ and $f \in G_a$, then $g^f \in N_a$.
To see this, choose $b \in P$ such that $b \not \in \cl(a g(a) f(g(a)))$
and 
$f(b) \not \in \cl(a g(a) f(g(a)))$. Let $h\in I$ such that $h(a)=b$. 
Since conjugation is a permutation of $I \setminus \{0 \}$,
and $X \cup f^{-1}(X)$ is finite, for each finite subset $X$ of $P$,
it suffices to show that $g^fh^f \in I$.
But this is clear since $gh \in I$.

It is easy to see that the action of $G_a$ is transitive, and even
sharply transitive by 2-determinedness.
Using 2-determinedness again, one shows that for each $g \in G$
there is $f \in N_a$ and $h \in G_a$ such that $g = fh$.
Since also $G_a \cap N_a = {0}$, we have that
$G = G_a \ltimes N_a$.

If $G_a$ is abelian, we define the structure of a field 
on $N_a$ as follows:
We let $N_a$ be the additive group of the field,
{\em i.e.} 
the addition $\oplus$ on $N_a$ is simply the group operation
of $N_a$ and $0$ its identity element. 
Now fix
an arbitrary element in 
$N_a \setminus \{ 0 \}$, which we denote by $1$ and
which will play the role of the identity.
For each $g \in N_a \setminus {0}$, let $f_g \in G_a$ be the
unique element such that $1^{f_g} = g$.
We define the multiplication $\otimes$ of elements $g,h \in N_a$ as follows:
$g \otimes h = h^{f_g}$.  It is easy to see that this makes $N_a$
into a field $K$.
This field carries an $\omega$-homogeneous pregeometry, and hence
it is algebraically closed by Theorem~\ref{t:ma}.

Now for the last sentence, let $b_i \in P$ for $i < 4$ be {\em distinct}
elements. We must show $g \in G$ such that $h(b_0) = b_1$ and $h(b_2) = b_3$.
Let $b_i' \in K (=N_a)$ such that $b_i'(a) = b_i$.
Then there are $f, g \in K$ such that
$f\cdot b_0' + g = b_1'$ and $f \cdot b_2' + g = b_3'$.
Let $f' \in G_a$ such that $1^{f'}(a) = f(a)$.
Then $g(f'(b_0)) = b_1$ and $g(f'(b_2)) = b_3$.
Hence, for all $b \in P_a$, we have $\cl(b) \setminus \cl(a) = \{b \}$
by $\omega$-homogeneity. 
Thus $P$ is a geometry and the action of $G$ on $P$ (as a set)
is sharply $2$-transitive. 
\end{proof}

We can obtain a geometry $P'$ from a pregeometry by taking
the quotient with the equivalence relation $E(x,y)$ given
by $\cl(x) = \cl(y)$, for $x, y \in P$.

\begin{proposition}~\label{p:regular}
Assume that $G$ $2$-acts sharply on the {\em geometry} $P$.
Then $G_a$ acts regularly on $(P_a)'$ and 
$N_a$ acts regularly on $P$.
\end{proposition}
\begin{proof}
The fact that $G_a$ acts regularly on $(P_a)'$ follows
from the last sentence of the previous proposition.
We now show that $N_a$ acts transitively on $P'$.
Suppose first that for some
$x \in P' \setminus \{ a \}$ the subgroup $Stab(x)$ of $N_a$ has
bounded index.  Then, since $N_a$ is connected (as it carries
an $\omega$-homogeneous pregeometry), we have $Stab(x) = N_a$,
and so $N_a x = \{ x \}$.
Let $y \in P \setminus \{ a \}$.  
$G_a$ acts transitively on $P_a$, so there is
$g \in G_a$ such that $gx = y$.
Then $N_a y = N_a gy = gN_a x = gx = y$, since $G_a$ normalises $N_a$.
But the action of $G$ on $P$ is $2$-determined, 
so the action of $N_a$ on $P$ is $2$-determined
and hence $N_a = \{ 0 \}$, a contradiction.
So, for each $x \in P \setminus \{ a \}$, 
the stabiliser $Stab(x)$ is proper.
An easy generalisation of Lemma~\ref{l:conn} 
therefore shows that it is finite-dimensional
(with respect to $\cl'$).
Since this holds for every $x \in P \setminus 
\{ a \}$, there is exactly one orbit
and $N_a$ acts transitively on $P \setminus \{ a \}$.
But $N_a a \not = \{ a \}$ since $G_a \cap N_a = \{ 0 \}$.
This implies that $N_a$ acts transitively on $P'$.

Now to see that the action of $N_a$ on $P$ is sharp, 
suppose that $gx= x$ for some $x \in P_a$.
Let $y \in P_a \setminus \cl(ax)$.
By transitivity, there is $h \in N_a$ such that $hx = y$.
Then $g y = ghx= hgx= hx = y$, since $N_a$ is abelian.
It follows that $g = 0$ by $2$-determinedness, so the action is regular.
\end{proof}

We can now obtain the full picture for groups acting on {\em geometries}.

\begin{theorem}~\label{t:big}
Let $G$ be a group $n$-acting on a geometry $P$.
Assume that $G$ admits hereditarily unique generics with respect to $\Sigma$.
Then, either there is an unbounded
non-classical $A$-invariant
subgroup of $G$, or $n = 1,2,3$ and 
\begin{enumerate}
\item
If $n = 1$, then $G$ is abelian and acts regularly
on $P$.
\item
If $n = 2$, then $P$ can be given the $A$-invariant structure
of an algebraically closed field $K$ (for $A \subseteq P$ finite), and 
the action of $G$ on $P$ is isomorphic to 
the affine action of $K^* \ltimes K^+$ on $K$.
\item
If $n = 3$, then 
$P \setminus \{ \infinity \}$ can be given the $A$-invariant 
structure
of an algebraically closed field $K$ (for some $\infinity \in P$
and $A \subseteq P$ finite), 
and the action of $G$ on $P$ is isomorphic to
the action of $\PGL_2(K)$ on the projective line $\mathbb{P}^1(K)$. 
\end{enumerate}
\end{theorem}
\begin{proof}
Suppose that there are no $A$-invariant unbounded non-classical 
subgroup of $P$,
for some finite $A$.
Then $(G_a)^0$ must be abelian, so that the action of $G$ on $P$ is
$n$-determined, by Fact~\ref{f:ndet}.
Thus $n = 1,2,3$ by Fact~\ref{f:3}.

For $n = 1$, the $G$ acts regularly on $P$, and hence carries a pregeometry
and must therefore be abelian (otherwise it is nonclassical).

For $n = 2$, notice that since $N_a$ acts regularly on $P$, we 
can endow $P$ with the algebraically closed field structure of $N_a$
by Proposition~\ref{t:field}.
The conclusion follows immediately. 

For $n = 3$, we follow \cite{Bu}, where some of this is done in the strongly
minimal case.
Choose a point $b \in P$ and call it $\infinity$.
Then by Proposition~\ref{t:field} the
$G_\infinity$ acts sharply $2$-transitively on the set $P_\infinity$, which is 
also a geometry, {\em i.e.} for each $b \in P_\infinity$, 
$\cl(b, \infinity) = \{ b \}$.
Hence by $\omega$-homogeneity of $P$, we have that 
$\cl(b, c) = \{ b, c \}$ for any $b, c \in P$, from which it
follows that the action of $G$ on the set $P$ is sharply $3$-transitive.

By (2) we can endow $P_\infinity$
with the structure of an algebraically closed field $K$.
Denote by $0$ and $1$ the identity elements of $K$.
Then $\{ \infinity, 0 , 1 \}$ is a set of dimension $3$.

Consider
$G_{\infinity, 0}$ which consists of those elements fixing
both $\infinity$ and $0$.
Then $G_{\infinity, 0 }$ carries an $\omega$-homogeneous pregeometry.
It is isomorphic to the multiplicative group
$K^*$.

Now let $\alpha$ be the unique element of $G$ sending
$(0,1, \infinity)$ to $(\infinity, 1, 0)$, which exists
since the action of $G$ on $P$ is sharply $3$-transitive.
Notice that $\alpha^2 = 1$.

We leave it to the reader to check that
conjugation by $\alpha$ induces an idempotent
automorphism $\sigma$ of $G_{\infinity, 0}$, which is not the identity.
Furthermore, $\sigma g = g^{-1}$ for each $g \in G_{\infinity, 0}$:
To see this, consider the proper
definable subgroup $B = \{ a \in G_{\infinity, 0} : \sigma(a) = a \}$
of $G_{\infinity, 0}$.  Then $B$ is $0$-dimensional in the pregeometry
$\cl'$ of $G_{\infinity, 0}$.  
Consider also $C = \{a \in G_{\infinity, 0} : \sigma(a) = a^{-1} \}$. 
Let $\tau : G_{\infinity, 0} \rightarrow G_{\infinity, 0}$ be
the homomorphism defined by $\tau(x) = \sigma(x) x^{-1}$.
Then for $x \in G_{\infinity, 0}$ we have
\[
\sigma(\tau(x))= \sigma^2(x)\sigma(x^{-1})=
x \sigma(x)^{-1} = \tau(x)^{-1},
\] 
so $\tau$ maps $G_{\infinity, 0}$ into $C$.
If $\tau(x) = \tau(y)$, then $x \in yB$, so $x \in \cl'(y)$
(in the pregeometry of $G_{\infinity, 0 }$).
It follows that the kernel of $\tau$ is finite dimensional,
and therefore $C = G_{\infinity, 0}$ (using essentially Lemma~\ref{l:ker}).

We can now complete the proof:
Given $x \in K^*$, choose $h \in G_{\infinity, 0}$ such that
$h1= x$.
Then $\alpha x = \alpha h 1 = h^{-1} \alpha 1 = h^{-1} 1 = x^{-1}$.
So $\alpha$ acts like an inversion on $K^*$. 
It follows that $G$ contains the group of automorphisms of $\mathbb{P}^1(K)$
generated by the affine transformations and inversion.
Hence $\PGL_2(K)$ embeds in $G$.
Since the action of $\PGL_2(K)$ and $G$ are both sharp $3$-actions, 
the embedding is all of $G$.

To see that $N_{0, \infinity}$ now carries the field $K$ and that the
action is as desired, it is enough to check that the correspondence
\[
N_{0,\infinity} \leftrightarrow G_{0, \infinity} \leftrightarrow 
P_{\infinity}
\]
commutes.
This follows from the following computation: For $0, 1, x \in P$, 
$1' \in N_{0, \infinity}$ chosen so that $1'(0) = 1$, and 
$h \in G_{0,\infinity}$ such that $h1= x$, we 
have
\[
h 1' h^{-1} (0) = h 1' 0 = h 1 = x.
\] 
Finally, going back from $P_\infinity$ to $P$, 
one checks easily that the action of $G$ on $P$ is isomorphic to
the action of $\PGL_2(K)$ on the projective line $\mathbb{P}^1(K)$. 
\end{proof}

\section{The stable homogeneous case}
We remind the reader of a few basic facts in homogeneous model
theory, which can be found in \cite{Sh:3}, \cite{HySh}, or \cite{GrLe}.
Let $L$ be a language and let $\bar{\kappa}$ be a suitably big cardinal.
Let $\mathfrak{C}$ be a {\em strongly $\bar{\kappa}$-homogeneous}
model, {\em i.e.} any elementary map 
$f : \mathfrak{C} \rightarrow \mathfrak{C}$ of size less than $\bar{\kappa}$
extends to an automorphism of $\mathfrak{C}$.
We denote by $\aut_A(\mathfrak{C})$ or $\aut(\mathfrak{C}/A)$
the group of automorphisms of $\mathfrak{C}$ fixing $A$ pointwise.
A set $Z$ will be called {\em $A$-invariant} if $Z$ is fixed setwise
by any automorphism $\sigma \in \aut(\mathfrak{C}/A)$.
This will be our substitute for definability; by homogeneity
of $\mathfrak{C}$ an $A$-invariant set is the 
disjunction of complete types over $A$. 

Let $D$ be the {\em diagram} of $\mathfrak{C}$, {\em i.e.} the 
set of complete $L$-types over the empty set realised by finite sequences
from $\mathfrak{C}$.
For $A \subseteq \mathfrak{C}$ we denote by
\[
S_D(A) = \{ p \in S(A) : \text{ For any $c \models p$ and $a \in A$
the type $\tp(ac/\emptyset) \in D$} \}.
\]
The homogeneity of $\mathfrak{C}$ has the following important consequence.
Let $p \in S(A)$ for $A \subseteq \mathfrak{C}$ with $|A| < \bar{\kappa}$.
The following conditions are equivalent:
\begin{itemize}
\item
$p \in S_D(A)$;
\item
$p$ is realised in $\mathfrak{C}$;
\item
$p \restriction B$ is realised in $\mathfrak{C}$ for each
finite $B \subseteq \mathfrak{C}$.
\end{itemize}
The equivalence of the second and third item is sometimes called
{\em weak compactness}, it is the chief reason why homogeneous model
theory is so well-behaved.

We will use $\mathfrak{C}$ as a universal domain; each set and model
will be assumed to be inside $\mathfrak{C}$ of size less than $\bar{\kappa}$,
satisfaction is taken with respect to $\mathfrak{C}$.
We will use the term {\em bounded} to mean `of size less than $\bar{\kappa}$'
and {\em unbounded} otherwise.  By abuse of language, a type is bounded
if its set of realisations is bounded.

We will work in the {\em stable} context.  We say 
that $\mathfrak{C}$ (or $D$) is {\em stable} if one of the following
equivalent conditions are satisfied:
\begin{fact}[Shelah] The following conditions are equivalent:
\begin{enumerate}
\item
For some cardinal $\lambda$, $D$ is {\em $\lambda$-stable},
{\em i.e.} $|S_D(A)| \leq \lambda$ for each $A \subseteq \mathfrak{C}$
of size $\lambda$.
\item
$D$ does not have the order property, {\em i.e.} there does not exist
a formula $\phi(x,y)$ such that for arbitrarily large $\lambda$ we have
$\{ a_i : i < \lambda \} \subseteq \mathfrak{C}$ such that
\[
\mathfrak{C} \models \phi(a_i, a_j) 
\quad
\text{if and only if}
\quad
i < j < \lambda.
\]
\item
There exists a cardinal $\kappa$ such that for each $p \in S_D(A)$
the type $p$ does not split over a subset $B \subseteq A$ of size
less than $\kappa$.
\end{enumerate}
\end{fact}
Recall that $p$ {\em splits over $B$} if there is $\phi(x,y) \in L$
and $c,d \in A$ with $\tp(c/B) = \tp(d/B)$ such that 
$\phi(x,c) \in p$ but $\neg \phi(x,d) \in p$.

Note that a diagram $D$ may be stable while the first order theory
of $\mathfrak{C}$ is unstable.  Further, in (3) the cardinal $\kappa$
is bounded by the first stability cardinal,
itself at most $\beth_{(2^{|L|})^+}$.

Nonsplitting provides a rudimentary independence relation in the context
of stable homogeneous model theory, but we will work primarily inside
the set of realisations of a quasiminimal type, where the independence
relation has a simpler form.
Recall that a type $p \in S_D(A)$ is {\em quasiminimal} (also
called {\em strongly minimal}) if it is unbounded but has a unique
unbounded (hence quasiminimal) extension to any $S_D(B)$, for $A \subseteq B$.
Quasiminimal types carry a pregeometry:

\begin{fact}~\label{f:preg}
Let $p$ be quasiminimal and let $P = p(\mathfrak{C})$.
Then $(P, \bcl)$, where for $a, B \subseteq P$ 
\[a \in \bcl_A(B)
\quad
\text{if} 
\quad 
\tp(a/A \cup B)
\text{ is bounded},
\] 
satisfies the axioms
of a pregeometry.
\end{fact}

We can therefore define $\dim(X/B)$ for $X \subseteq P = p(\mathfrak{C})$
and $B \subseteq \mathfrak{C}$.
This induces a dependence relation $\nonfork_{}$ as follows:
\[
a \nonfork_B C,
\]
for $a \in P$ a finite sequence, and $B, C \subseteq \mathfrak{C}$ if
and only if 
\[
\dim(a/B) = \dim(a/B \cup C).
\]
We write $\fork_{}$ for the negation of $\nonfork_{}$.
The following lemma follows easily.

\begin{lemma}  Let $a, b \in P$ be finite sequences, and
$B \subseteq C \subseteq D \subseteq E \subseteq \mathfrak{C}$.
\begin{enumerate}
\item (Finite Character)
If $a \fork_B C$, then there exists a finite $C' \subseteq C$
such that $a \fork_B C'$. 
\item (Monotonicity)
If $a \nonfork_B E$ then $a \nonfork_C D$.
\item (Transitivity)
$a \nonfork_B D$ and $a \nonfork_D E$ if and only if $a \nonfork_B E$.
\item (Symmetry)
$a \nonfork_B b$ if and only if $b \nonfork_B a$.
\end{enumerate}
\end{lemma} 

This dependence relation (though defined only some sets in $\mathfrak{C}$)
allows us to extend much of the theory of forking.

\relax From now until Theorem~\ref{t:main},
we make the following hypothesis:
\begin{hypothesis}\label{h:basic}
Let $\mathfrak{C}$ be stable.
Let $p$, $q \in S_D(A)$ be unbounded, with $p$ quasiminimal.
Let $n < \omega$ be such that:
\begin{enumerate}
\item
For any independent sequence $(a_0, \dots, a_{n-1})$
of realisations of $p$
and any (finite) set $C$ of realisations of $q$ we
have
\[
\dim(a_0, \dots, a_{n-1}/A) = \dim(a_0, \dots, a_{n-1}/ A \cup C).
\]
\item
For some independent sequence $(a_0, \dots, a_n)$ of realisations
of $p$ there is a finite set $C$ of realisations of $q$ such that
\[
\dim(a_0, \dots, a_n/A) > \dim(a_0, \dots, a_n/ A \cup C).
\]
\end{enumerate}
\end{hypothesis}

\begin{remark}
In case we are in the $\omega$-stable~\cite{Le:1} or even 
the superstable~\cite{HyLe} case, there is a dependence relation
on all the subsets, induced by a rank, which satisfies  
many of the properties of forking (symmetry and extension
only over certain sets, however).  This dependence relation,
which coincides with the one defined when both make sense, 
allows us to develop orthogonality calculus in much the
same way as the first order setting, and would have enabled us
to phrase the conditions (1) and (2) in the same way as the one we
phrased for Hrushovski's theorem.
Without canonical bases, however, it is not clear that the,
apparently weaker, condition that $p^n$ is weakly orthogonal to 
$q^\omega$ implies (1).
\end{remark}

We now make the pregeometry $P$ into a {\em geometry} $P/E$ by 
considering the equivalence relation $E$ on elements of $P$ given
by 
\[
E(x,y)
\quad
\text{if and only if} 
\quad
\bcl_A(x) = \bcl_A(y).
\]

We now proceed with the construction.
Before we start, recall that the notion of interpretation we
use in this context is like the first order notion, except
that we replace definable sets by invariant sets (see Definition~\ref{d:int}).

Let $Q = q(\mathfrak{C})$.
The group we are going to interpret is the following:
\[
\aut_{Q \cup A}(P/E).
\]
The group $\aut_{Q \cup A}(P/E)$ is the group of permutations
of the geometry obtained from $P$, which are induced by automorphisms
of $\mathfrak{C}$ fixing $Q \cup A$ pointwise.  
There is a natural action of this group on the geometry $P/E$.
We will show in this section that the action has rank $n$,
is $(n+1)$-determined.  Furthermore, considering the automorphisms
induced from $\aut_{A}(\mathfrak{C})$, we have a group acting on a geometry
in the sense of the previous section.
By restricting the group of automorphism to those induced
by the group of strong automorphisms
$\Saut_{A}(\mathfrak{C})$,
we will show in addition that this group is stationary
and admits hereditarily unique generics.
The conclusion will then follow easily from the last theorem
of the previous section.

We now give the construction more precisely.
\begin{notation}
We denote by $\aut(P/A \cup Q)$ the group
of permutations of $P$ which extend to an automorphism
of $\mathfrak{C}$ fixing $A \cup Q$.
\end{notation}

Then $\aut(P/A\cup Q)$
acts on $P$ in the natural way.  Moreover, each $\sigma \in \aut(P/A\cup Q)$ 
induces a unique permutation on $P/E$, which we denote by $\sigma/E$.
We now define the group that we will interpret:
\begin{definition}
Let $G$ be the group consisting of the
permutations $\sigma/E$
of $P/E$ induced by elements $\sigma \in \aut(P/A \cup Q)$.
\end{definition}

Since $Q$ is unbounded, $\aut_{A \cup Q}(\mathfrak{C})$
could be trivial
(this is the case even in the first order case if the theory is
not stable).
The next lemma shows that this is not the case under stability
of $\mathfrak{C}$.
By abuse of notation, we write 
\[
\tp(a/A \cup Q) = \tp(b/A \cup Q),
\] 
if
$\tp(a/AC) = \tp(b/AC)$ for any bounded $C \subseteq Q$.

\begin{lemma}\label{l:aut} Let $a, b$ be bounded sequences in $\mathfrak{C}$
such that 
\[
\tp(a/A \cup Q) = \tp(b/A \cup Q).
\]
Then there exists $\sigma \in \aut(\mathfrak{C})$ sending $a$ to $b$
which is the identity
on $A \cup Q$.
\end{lemma}
\begin{proof}
By induction, it is enough to prove that 
for all $a' \in \mathfrak{C}$, there is $b' \in \mathfrak{C}$
such that $\tp(aa'/A \cup Q) = \tp(bb'/A \cup Q)$.

Let $a' \in \mathfrak{C}$.
We claim that there exists a bounded 
$B \subseteq Q$ such that
for all $C \subseteq Q$ bounded, we have 
$\tp(aa'/ABC)$ does not split over $AB$.

Otherwise, for any $\lambda$, 
we can inductively construct an increasing sequence
of bounded sets $(C_i : i < \lambda)$ such
that $\tp(aa'/C_{i+1})$ does not split over $C_i$.
This contradicts stability (such a chain must stop at the first
stability cardinal).

Now let $\sigma \in \aut_{A \cup B}(\mathfrak{C})$ sending $a$ to $b$
and let $b' = \sigma(a')$.
We claim that $\tp(aa'/A \cup Q) = \tp(bb'/A \cup Q)$.
If not, let $\phi(x,y,c) \in \tp(aa'/A \cup Q)$ and
$\neg \phi(x,y,c) \in \tp(bb'/A \cup Q)$.
Then, $\phi(x,y,c), \neg \phi(x,y, \sigma(c)) \in \tp(aa'/A \cup Bc\sigma(c))$,
and therefore $\tp(a a'/AB c \sigma(c))$ splits over $AB$, a contradiction.
\end{proof}

It follows that the action of $\aut(P/A \cup Q)$ on $P$, and a fortiori 
the action of $G$ on $P/E$, 
has some transitivity properties.
The next corollary implies that the action of $G$ on $P$ has
rank $n$ (condition (2) in Hypothesis~\ref{h:basic} prevents
two distinct independent $(n+1)$-tuples of realisations of $p$
from
being automorphic over $A \cup Q$).

\begin{corollary}\label{c:aut} 
For any independent
$a, b \in P^n$, there is $g \in G$
such that $g(a/E) = b/E$.
\end{corollary}
\begin{proof}
By assumption (1) $\dim(a/A \cup C) = \dim(b/ A \cup C) = n$, for each
finite $C \subseteq Q$.
By uniqueness of unbounded extensions, we have that 
$\tp(a/A \cup C) = \tp(b/A \cup C)$ for each finite $C \subseteq Q$.
It follows that $\tp(a/A \cup Q) = \tp(b/A \cup Q)$ so by the previous
lemma, there is $\sigma \in \aut(\mathfrak{C}/A \cup C)$ such
that $\sigma(a) = b$.
Then $g = \sigma/E$.
\end{proof}

The next few lemmas are in preparation to show that the 
action is $(n+1)$-determined.
We first give a condition ensuring that two elements of $G$ coincide.

\begin{lemma}\label{l:1.4} Let $\sigma, \tau \in \aut(P/A \cup Q)$.
Let $a_i, b_i \in P$, for $i < 2n$ be such that 
\[
\sigma(a_i) = b_i =
\tau(a_i),
\quad
\text{for $i < 2n$}.
\]
Assume further that
\[
a_i \nonfork_A \{ a_j, b_j : j < i \}
\quad
\text{and}
\quad
b_i \nonfork_A \{ a_j, b_j : j < i \},
\quad
\text{for $i < 2n$}.
\]
Let $c \in P$ be such that $c, \sigma(c), \tau(c) \not \in 
\bcl_A(\{ a_i, b_i : i < 2n \})$.
Then 
\[
\sigma(c)/E) = \tau(c)/E).
\]
\end{lemma}

\begin{proof}
Let
$\bar{a} = (a_i : i < 2n)$ and $\bar{b}= (b_i : i < 2n)$
satisfy the independence requirement, and $\sigma(a_i) = b_i = \tau(a_i)$, for
$i < 2n$.
Assume, for a contradiction, that $c \in P$ is as above but 
$\sigma(c)/E \not = \tau(c)/E$.

We now establish a few properties:
\begin{enumerate}
\item
$\sigma(c) \fork_{A \cup \{ a_i, b_i : i < n \}} c$,
\item
$\sigma(c) \fork_{A \cup \{ a_i, b_i : n \leq i < 2n \}} c$,
\item
$\tau(c) \fork_{A \cup \{ a_i, b_i : i < n \}} c$,
\item
$\tau(c) \fork_{A\cup \{ a_i, b_i : n \leq i < 2n \}} c$ 
\end{enumerate}
All these statements are proved the same way, so we only show (1):
Suppose, for a contradiction, 
that $\sigma(c) \nonfork_{A \cup \{ a_i, b_i : i < n \}} c$.
By Hypothesis~\ref{h:basic}, there is a finite $C \subseteq Q$
such that $\dim(c a_0 \dots a_{n-1} /A \cup C) \leq n$,
{\em i.e.}
\begin{equation} \tag{*}
c a_0 \dots a_{n-1}  \fork_A C.
\end{equation}
Let $c' \in P$ be such that 
$c' \not \in \bcl_A(C \cup c \cup \{ a_i, b_i : i < n\})$.
By assumption, $\sigma(c) \not \in \bcl_A(c \cup \{ a_i, b_i : i < n\})$.
Hence, there exists an automorphism $f$ of
$\mathfrak{C}$ such that $f(c') = \sigma(c)$ 
which is the identity on $A \cup
c \cup
\{ a_i, b_i : i < n \}$.
Then by using $f$ on (*), we obtain 
\[
c a_0 \dots a_{n-1}  \fork_A f(C)  
\]
On the other hand, 
$\sigma(c) \not \in \bcl_A(f(C) \cup c \cup \{ a_i, b_i : i < n \})$,
since $\sigma(c) = f(c')$.
But by Hypothesis~\ref{h:basic} we have 
\[
b_0 \dots b_{n-1} \nonfork_A f(C).
\]
Together these imply
\[
\sigma(c) b_0 \dots b_{n-1} \nonfork_A f(C). 
\]
But this contradicts (*), since $\sigma$ fixes $f(C) \subseteq Q$.

We now prove another set of properties:
\begin{enumerate}
\item[(5)]
$\sigma(c) \fork_{A \cup \{ b_i : i < n \}} \tau(c)$
\item[(6)]
$\sigma(c) \fork_{A \cup \{ b_i : n \leq i < 2n \}} \tau(c)$
\end{enumerate}
These are again proved similarly using the fact that for all finite 
$C \subseteq Q$
\[
\sigma(c) \cup \{ b_i : i < n \} \fork_A C
\quad
\text{ if and only if }
\quad 
\tau(c) \cup \{ b_i : i < n \} \fork_A C.
\]

We can now finish the claim:  By (6) we have that
\[
\{ b_i : n \leq i < 2n \} \fork_A \sigma(c) \tau(c).
\]
This implies that 
\[
\{ b_i : n \leq i < 2n \} \fork_{A \cup \{ b_i : i
< n \}} \sigma(c) \tau(c),
\]
since $\{ b_i : i < 2n \}$ are independent.
By (5) using the fact that $(P, \bcl_A)$ is a pregeometry, we therefore
derive that
\[
\{ b_i : n \leq i < 2n \} \fork_{A \cup \{ b_i : i
< n \}} \sigma(c).
\]
But this contradicts the fact that $\sigma(c) \nonfork_A \bar{a} \bar{b}$.
\end{proof}

\begin{lemma}\label{l:2}
Let $\sigma \in \aut(P/A \cup Q)$ and $\bar{a} \in P^{n+1}$ be independent.
If 
\[
\sigma(a_i)/E = a_i/E,
\quad
\text{for each $i \leq n$},
\] 
and $c \not \in \bcl_A(\bar{a}\sigma(\bar{a}))$,
then 
\[
\sigma(c)/E = c/E.
\]
\end{lemma}

\begin{proof}  
Suppose, for a contradiction, that the conclusion fails.
Let $c \in P$, with $c \not \in \bcl(A \bar{a} \sigma(\bar{a})$ 
such that $\sigma(c) \not 
\in \bcl_A(c)$.
Choose $a_i$, for $n < i < 2n+1$ such that
$a_i$ and $\sigma(a_i)$ satisfy the assumptions of Lemma~\ref{l:1.4}
and 
\[
c a_0 \nonfork_{A} \{ a_i, \sigma(a_i) : 0 < i < 2n+1 \}. 
\]
This is possible: To see this, assume that
we have found $a_j$ and $\sigma(a_j)$, for $j < i$, satisfying the requirement.
For each $k \leq j$, choose $a'_k$ such that 
\[
a'_k \nonfork_A \{ a'_\ell : \ell < k \} \cup \{ a_j , \sigma(a_j) : j < i \}.
\]
Then,
\[
\dim(\{ a_j : j < i \} \cup \{ a'_k : k \leq i \}) = 2i+1.
\]
Since $\sigma$ extends to an automorphism of $\mathfrak{C}$, we must have also
\[
\dim(\{ \sigma(a_j) : j < i \} \cup \{ \sigma(a'_k) : k \leq i \}) = 2i+1.
\]
Hence, for some $k \leq i$ we have
\[
\sigma(a'_k) \nonfork_A \{ a_j , \sigma(a_j) : j < i \},
\]
and we can let $a_i = a'_k$ and $\sigma(a_i) = \sigma(a'_k)$.

Then there is an automorphism $f$ of $\mathfrak{C}$ which sends $c$ to $a_0$
and is the
identity on $A \cup \{ a_i, \sigma(a_i) : 0 < i < 2n+1\}$.
Then $\sigma$ and $f^{-1} \circ \sigma \circ f$ contradict Lemma~\ref{l:1.4}.
\end{proof}

We can now obtain:

\begin{lemma}\label{l:3}
Let $\sigma \in \aut(P/A \cup Q)$.
Assume that $(a_i)_{i \leq n} \in P^{n+1}$ is independent
and 
$\sigma(a_i)/E = a_i/E$, for $i \leq n$.
Then $\sigma /E$ is the identity in $G$.
\end{lemma}
\begin{proof} Let
$\bar{a} \in P^{n+1}$ be independent.
Choose $a \in P$ arbitrary and $a_i' \in P$ for $i < n+1$
such that $a'_i \not \in \bcl(Aa_0, \dots, a_n, a_0', \dots a_i')$
for each $i < n+1$.
By the previous lemma, $\sigma(a_i') \in \bcl_A(a_i')$ for each $i < n+1$.
Hence $\sigma(a) \in \bcl_A(a)$ by another application of the lemma.
\end{proof}

The next corollary follows by applying the lemma to $\tau^{-1} \circ \sigma$.
Together with Corollary~\ref{c:aut}, it shows that the action
of $G$ on $P/E$ is an $n$-action.

\begin{corollary}
Let $\sigma, \tau \in \aut(P/A \cup Q)$ and assume there is an 
$(n+1)$-dimensional
subset $X$ of $P/E$ on which $\sigma/E$ and $\tau/E$ agree.
Then $\sigma/E = \tau/E$.
\end{corollary}

We now consider automorphisms of this group action.
Let $\sigma \in \aut_A(\mathfrak{C})$.
Then, $f$ induces 
an automorphism $\sigma'$ of the group action as follows:
$\sigma'$ is $\sigma/E$ on $P/E$, and for $g \in G$
we let $\sigma'(g)(a/E)= 
\sigma(\tau(\sigma^{-1}(a)))/E$, where $\tau$ is such that $\tau/E = g$.
It is easy to verify that 
\[
\sigma' : G \rightarrow G
\]
is an automorphism
of $G$ (as $\sigma \circ \tau \circ \sigma^{-1} 
\in \aut_{Q \cup A}(\mathfrak{C})$
if $\tau \in \aut_{Q \cup A}(\mathfrak{C})$, 
and both $P$ and $Q$ are $A$-invariant).
Finally, one checks directly that $\sigma'$ preserves the action.

For stationarity, it is more convenient to 
consider strong automorphisms.
Recall that two sequences $a, b \in \mathfrak{C}$ have the same {\em Lascar
strong types} over $C$, written $\Lstp(a/C) = \Lstp(b/C)$,
if $E(a,b)$ holds for any $C$-invariant equivalence relation $E$ with
a bounded number of classes.
An automorphism $f \in Aut(\mathfrak{C}/C)$ is called {\em strong}
if $\Lstp(a/C) = \Lstp(f(a)/C)$ for any $a \in \mathfrak{C}$.
We denote by $\Saut(\mathfrak{C}/C)$ or $\Saut_C(\mathfrak{C})$
the group of strong automorphisms fixing $C$ pointwise.
We let $\Sigma = \{ \sigma' : \sigma \in \Saut_A(\mathfrak{C}) \}$.
The reader is referred to \cite{HySh} or \cite{BuLe} for more
details.

First, we show that the action is $\omega$-homogeneous 
with respect to $\Sigma$.

\begin{lemma} If $X \subseteq P/E$ is 
finite and $x, y \in P/E$ are outside $\bcl_A(X)$, then there is 
an automorphism $\sigma \in \Sigma$ of the group action sending
$x$ to $y$ which is the identity on $X$.
\end{lemma}  
\begin{proof}
By uniqueness of unbounded extensions, there is an automorphism $\sigma \in 
\Saut(\mathfrak{C})$
fixing $A \cup X$ pointwise and sending $x$ to $y$.
The automorphism $\sigma'$ is as desired.
\end{proof}

We are now able to show the stationarity of $G$.

\begin{proposition}~\label{p:stat} $G$ is stationary 
with respect to $\Sigma$. 
\end{proposition}
\begin{proof}
First, notice that the number of strong types is bounded by stability.
Now, let $g \in G$ be generic over the bounded set $X$ and let 
$\bar{x} \in  P^n$
be an independent sequence witnessing this, {\em i.e.}
\[
\dim(\bar{x} g(\bar{x})/X) =  2n.
\]
If $x' \in P$ is such that $\dim(\bar{x} x' /X) = n+1$, then 
$\dim(\bar{x} x'g(\bar{x})g(x)'/X) =  2n +1$.
By quasiminimality of $p$, this implies that
\[
\bar{x} x'g(\bar{x})g(x') \nonfork_A X.
\] 
Now let $h \in G$ be also generic over $X$ and such that
$\sigma(g) = h$ with $\sigma \in \Sigma$.
For $\bar{y}, y'$ witnessing the genericity of $h$ as above,
we have
\[
\bar{y} y'h(\bar{y})h(y') \nonfork_A X.
\]
Hence, by stationarity of Lascar strong types we
have $\Lstp(\bar{x} x'g(\bar{x})g(x')/A X) = 
\Lstp(\bar{y} y'h(\bar{y})h(y')/AX)$.
Thus, 
there is $\tau$, a strong automorphism of $\mathfrak{C}$ fixing $A \cup X$
pointwise, such that
$\tau(\bar{x}x'g(\bar{x}g(x')) = \bar{y}y'h(\bar{y})h(y')$.
Then, $\tau'(g) = h$ ($\tau' \in \Sigma$) since the action is 
$(n+1)$-determined.
\end{proof}

The previous proposition implies that $G^0$ has unique generics,
but we can prove more:

\begin{proposition}\label{p:gen}
$G^0$ admits hereditarily unique generics with respect to $\Sigma$.
\end{proposition}
\begin{proof}
For any independent $k$-tuple $a \in P/E$ with $k <n$,
consider the $\Sigma_a$-homogeneous $(n-k)$-action $G_a$ on $P/E$.
Instead of $\Sigma_a$, consider the smaller group $\Sigma_a'$
consisting of $\sigma'$ for strong automorphisms of $\mathfrak{C}$
fixing $Aa$ and preserving strong types over $Aa$.  
Then, as in the proof of the previous proposition,
$G_a$ is stationary with respect to $\Sigma_a'$, which implies
that the connected component $G_a'$ of $G_a$ (defined with $\Sigma_a'$)
has unique generics with respect to restriction
of automorphisms in $\Sigma_a'$ by Theorem~\ref{t:unique}.
But, there are even more automorphisms in $\Sigma_a$ so $G_a'$
has unique generics with respect to restriction of automorphisms in
$\Sigma_a$.
By definition, this means that $G$ admits hereditarily unique generics.
\end{proof}

We now show that $G$ is interpretable in $\mathfrak{C}$.
We recall the definition of interpretable group in this context.

\begin{definition}\label{d:int}
A group $(G,\cdot)$ {\em interpretable in
$\mathfrak{C}$} if there is a (bounded) subset  $B \subseteq \mathfrak{C}$ 
and an unbounded set $U \subseteq \mathfrak{C}^k$
(for some $k < \omega$), an equivalence
relation $E$ on $U$, and a binary function $*$ on $U/E$
which are $B$-invariant and such that 
$(G,\cdot)$ is isomorphic to $(U/E,*)$.
\end{definition}

We can now prove:

\begin{proposition}\label{c:main} The group $G$ is interpretable
in $\mathfrak{C}$.
\end{proposition}
\begin{proof}
This follows from the $(n+1)$-determinacy of the group action.
Fix $a$ an independent $(n+1)$-tuple of elements of $P/E$.
Let $B = Aa$.

We let $U/E \subseteq P^{n+1}/E$ consist of those $b \in P^{n+1}/E$
such that $ga = b$ for some $G$.
Then, this set is $B$-invariant since if $b \in P^{n+1}/E$ and $\sigma \in
\aut_B(\mathfrak{C})$, then $\sigma'(g) \in G$ and sends $a$ to $\sigma(b)$
(recall that $\sigma'$
is the automorphism of the group action induced by $\sigma$).

We now define $b_1 * b_2 = b_3$ on $U/E$, if whenever
$g_\ell \in G$ such that $g_\ell(a) = b_\ell$, then
$g_1 \circ g_2 = g_3$.
This is well-defined by $(n+1)$-determinacy and the definition
of $U/E$.  Furthermore, the binary function $*$ is $B$-invariant.
It is clear that $(U/E,*)$ is isomorphic to $G$.
\end{proof}

\begin{remark}
As we pointed out, by homogeneity of $\mathfrak{C}$,
any $B$-invariant set is equivalent to a disjunction of 
complete types over $A$.
So, for example, if $B$ is finite, 
$E$ and $U$ is expressible by formulas in 
$L_{\lambda+, \omega}$, where $\lambda = |S^D(B)|$.
\end{remark}

It follows from the same proof that 
$G^0$ is interpretable in $\mathfrak{C}$,
and similarly $G_a$ and $(G_a)^0$ are interpretable for
any independent $k$-tuple $a$ in $P/E$ with $k <n$.

\begin{remark}
If we choose $p$ to be regular (with respect to, say, strong splitting),
we can still interpret a group $G$, exactly as we have in the case
of $p$ quasiminimal.
We have used the fact that the dependence relation is given by bounded
closure only to ensure the stationarity of $G$, and to obtain a field.
\end{remark}

We can now prove the main theorem.
We restate the hypotheses for completeness.

\begin{theorem}~\label{t:main}
Let $\mathfrak{C}$ be a large, homogeneous model of
a stable diagram $D$.
Let $p, q \in S_D(A)$ be unbounded with $p$ quasiminimal.
Assume that there is $n \in \omega$ such that 
\begin{enumerate}
\item
For any independent $n$-tuple $(a_0, \dots, a_{n-1})$
of realisations of $p$
and any finite set $C$ of realisations of $q$ we
have
\[
\dim(a_0, \dots, a_{n-1}/ A \cup C) = n.
\]
\item
For some independent sequence $(a_0, \dots, a_n)$ of realisations
of $p$ there is a finite set $C$ of realisations of $q$ such that
\[
\dim(a_0, \dots, a_n/ A \cup C) < n+1.
\]
\end{enumerate}
Then $\mathfrak{C}$ interprets a group $G$ which acts on the geometry
$P'$ obtained from $P$.
Furthermore, either $\mathfrak{C}$ interprets a non-classical group,
or $n \leq 3$ and 
\begin{itemize}
\item
If $n = 1$, then $G$ is abelian and acts
regularly on $P'$;
\item
If $n = 2$, the action of $G$ on $P'$ is isomorphic to the affine
action of $K^+ \rtimes K^*$ on the algebraically closed field $K$.
\item
If $n = 3$, the action of $G$ on $P'$ is isomorphic to the action
of $\PGL_2(K)$ on the projective line $\mathbb{P}^1(K)$ of the
algebraically closed field $K$.
\end{itemize}
\end{theorem}

\begin{proof}
The group $G$ is interpretable in $\mathfrak{C}$ by Proposition~\ref{c:main}.
This group acts on the geometry $P/E$; the action has rank $n$ and
is $(n+1)$-determined.
Furthermore, $G^0$ admits hereditarily unique generics
with respect to set of automorphisms $\Sigma$
induced by strong automorphisms of $\mathfrak{C}$.
Working now with the connected group $G^0$, which is invariant and therefore
interpretable,  
the conclusion follows from Theorem~\ref{t:big}.
\end{proof}

\begin{question}
The only point where we use quasiminimality is in showing that 
$G$ admits hereditarily unique generics.  Is it possible to do
this for regular types, say in the superstable case?
\end{question}

\section{The excellent case}

Here we consider a class $\mathcal{K}$ of atomic models of a countable
first order theory, {\em i.e.} $D$ is the set of isolated types over the
empty set.
We assume that $\K$ is {\em excellent} (see \cite{Sh:87a}, \cite{Sh:87b},
\cite{GrHa} or \cite{Le:3} for the basics of excellence).
We will use the notation $S_D(A)$ and splitting, which have been defined
in the previous section.

Excellence lives in the {\em $\omega$-stable} context,
{\em i.e.} $S_D(M)$ is countable, for any countable $M \in \K$.
This notion of $\omega$-stability is strictly weaker than the corresponding
notion given in the previous section; in the excellent, non-homogeneous
case, there are countable atomic sets $A$ such that $S_D(A)$ is uncountable. 
{\em Splitting} provides an dependence relation between sets,
which satisfies all the usual axioms of forking, provided we 
only work over models in $\K$.
For each $p \in S_D(M)$, for $M \in \K$, there is a finite
$B \subseteq M$ such that $p$ does not split over $B$.
Moreover, if $N \in \K$ extends $M$ then $p$ has a unique
extension in $S_D(N)$ which does not split over $B$.
Types with a unique nonsplitting extension are called {\em stationary}.

Excellence is a requirement on the existence of primary models,
{\em i.e.} a model $M \in \K$ is {\em primary over $A$},
if $M = A \cup \{ a_i : i < \lambda \}$ and for each $i < \lambda$
the type $\tp(a_i/A \cup \{ a_j : j < i \})$ is isolated.
Primary models are prime in $\K$.
The following fact is due to Shelah~\cite{Sh:87a}, \cite{Sh:87b}:
\begin{fact}[Shelah]
Assume that $\K$ is excellent.
\begin{enumerate}
\item
If $A$ is a finite atomic set, then there is a primary
model $M \in \K$ over $A$.
\item
If $M \in \K$ and $p \in S_D(M)$, then
for each $a \models p$, there is a primary
model over $M \cup a$.
\end{enumerate}
\end{fact}

We will use {\em full models} as universal domains
(in general $\K$ does not contain uncountable homogeneous
models).  
The existence of arbitrarily large full models follows from excellence.
They have the following properties (see again~\cite{Sh:87a} and~\cite{Sh:87a}):
\begin{fact}[Shelah] Let $M$ be a full model of uncountable size 
$\bar{\kappa}$.
\begin{enumerate}
\item
$M$ is $\omega$-homogeneous.
\item
$M$ is {\em model-homogeneous}, {\em i.e.} if
$a, b \in M$ have the same type over $N \prec M$ with $\| N \| < \bar{\kappa}$,
then there is an automorphism of $M$ fixing $N$ sending $a$ to $b$.
\item
$M$ realises any $p \in S_D(N)$ with $N \prec M$ of size less than 
$\bar{\kappa}$.
\end{enumerate}
\end{fact}

We work inside a full
$\mathfrak{C}$ of size $\bar{\kappa}$,
for some suitably big cardinal $\bar{\kappa}$.
All sets and models will be assumed to be inside $\mathfrak{C}$
of size less than $\bar{\kappa}$, unless otherwise specified.
The previous fact shows that all types over finite sets, and all
stationary types of size less than $\bar{\kappa}$ are realised
in $\mathfrak{C}$.

Since the automorphism group of $\mathfrak{C}$ is not as rich as in
the homogeneous case, it will be necessary to consider another closure
operator:
For all $X \subseteq \mathfrak{C}$ and $a \in M$, we define the
{\em essential closure} of $X$, written $\ecl(X)$ by 
\[
a \in \ecl(X),
\quad
\text{ if $a \in M$ for each $M \prec \mathfrak{C}$ containing $X$}.
\]

As usual, for $B \subseteq \mathfrak{C}$, 
we write $\ecl_B(X)$ for the closure operator on subsets $X$ of $\mathfrak{C}$
given by $\ecl(X \cup B)$.
Over finite sets, 
essential closure coincides with bounded closure, because
of the existence of primary models. 
Also, it is easy to check that $X \subseteq \ecl_B(X) = \ecl_B(\ecl_B(X))$,
for each $X, B \subseteq \mathfrak{C}$.
Furthermore, $X \subseteq Y$ implies that $\ecl_B(X) \subseteq \ecl_B(Y)$.

Again we consider a {\em quasiminimal} type $p \in S_D(A)$,
{\em i.e.} $p(\mathfrak{C})$ is unbounded and there is a unique
unbounded extension of $p$ over each subset of $\mathfrak{C}$.
Since the language is countable in this case, and we have $\omega$-stability,
the bounded closure of a countable set is countable.
Bounded closure satisfies exchange on the set of realisations
of $p$ (see \cite{Le:3}).  This holds also for essential closure.
\begin{lemma}  Let $p \in S_D(A)$ be quasiminimal.
Suppose that $a, b \models p$ are such that
$a \in \ecl_B(Xb) \setminus \ecl_B(X)$.
Then $b \in \ecl_B(Xa)$.
\end{lemma}
\begin{proof}
Suppose not, and let $M \prec \mathfrak{C}$ containing $A \cup B \cup X \cup a$
such that $b \not \in M$.  
Let $N$ containing $A \cup B \cup X$ such that $a \not \in N$.
In particular $a \not \in \bcl_B(N)$ and $a \in \ecl_B(Nb)$.
Let $b' \in \mathfrak{C}$ realise the unique free extension of $p$ over 
$M \cup N$.  Then $\tp(b/M) = \tp(b'/M)$ since there is a unique big extension
of $p$ over $M$.  It follows that there
exists $f \in \aut(\mathfrak{C}/M)$ such that $f(b) = b'$.
Let $N' = f(N)$.  Then $b' \not \in \bcl_B(N'a)$. 
On the other hand, we have $a \in \ecl_B(Nb) \setminus \ecl_B(N)$ 
by monotonicity
and choice of $N$, so $a \in \ecl_B(N'b') \setminus \ecl_B(N')$.
But, then $a \in \bcl_B(N'b') \setminus \bcl(N')$ 
(if $a \not \in \bcl_B(N'b')$,
then $a \not \in N'(b')$, for some (all) primary models over $N' \cup b'$).
But this is a contradiction.
\end{proof}

It follows from the previous lemma that the closure relation $\ecl_B$ 
satisfies the axioms of a pregeometry on the {\em finite} subsets of 
$P = p(\mathfrak{C})$, when $p$ is quasiminimal. 

Thus, for finite subsets $X \subseteq P$, and any set 
$B \subseteq \mathfrak{C}$, we can define $\dim(X/B)$ using
the closure operator $\ecl_B$.
We will now use the independence relation $\nonfork_{}$ as follows:
\[
a \nonfork_B C,
\]
for $a \in P$ a finite sequence, and $B, C \subseteq \mathfrak{C}$ if
and only if 
\[
\dim(a/B) = \dim(a/B \cup C).
\]

The following lemma follows easily. 

\begin{lemma}  Let $a, b \in P$ be finite sequences, and
$B \subseteq C \subseteq D \subseteq E \subseteq \mathfrak{C}$.
\begin{enumerate}
\item (Monotonicity)
If $a \nonfork_B E$ then $a \nonfork_C D$.
\item (Transitivity)
$a \nonfork_B D$ and $a \nonfork_D E$ if and only if $a \nonfork_B E$.
\item (Symmetry)
$a \nonfork_B b$ if and only if $b \nonfork_B a$.
\end{enumerate}
\end{lemma}

\relax From now until Theorem~\ref{t:main2},
we now make a hypothesis similar to Hypothesis~\ref{h:basic},
except that $A$ is chosen finite and the witness $C$ is allowed
to be countable (the reason is that we do not have
finite character in the right hand-side argument of $\nonfork_{}$).
Since we work over finite sets, notice that $p$ and $q$ below are
actually equivalent to formulas over $A$.

\begin{hypothesis}\label{h:hyp2}
Let $\mathfrak{C}$ be a large full model of an excellent class $\K$.
Let $A \subseteq \mathfrak{C}$ be finite.
Let $p$, $q \in S_D(A)$ be unbounded with $p$ quasiminimal.
Let $n < \omega$.
Assume that
\begin{enumerate}
\item
For any independent sequence $(a_0, \dots, a_{n-1})$
of realisations of $p$
and any countable set $C$ of realisations of $q$ we
have
\[
\dim(a_0, \dots, a_{n-1}/A) = \dim(a_0, \dots, a_{n-1}/ A \cup C).
\]
\item
For some independent sequence $(a_0, \dots, a_n)$ of realisations
of $p$ there is a countable set $C$ of realisations of $q$ such that
\[
\dim(a_0, \dots, a_n/A) > \dim(a_0, \dots, a_n/ A \cup C).
\]
\end{enumerate}
\end{hypothesis}

Write $P = p(\mathfrak{C})$ and $Q = q(\mathfrak{C})$, as in the previous
section.
Then, $P$ carries a pregeometry with respect to bounded closure,
which coincides with essential closure over finite sets.
Thus, when we speak about finite sets or sequences in $P$,
the term independent is unambiguous.
We make $P$ into a geometry $P/E$ by considering
the $A$-invariant
equivalence relation 
\[
E(x,y), 
\quad
\text{defined by}
\quad 
\bcl_A(x) = \bcl_A(y).
\]

The group we will interpret in this section is
defined slightly differently, because of the lack of homogeneity
(in the homogeneous case, they coincide).
We will consider the group $G$ of all permutations $g$ of $P/E$
with the property that for each countable $C \subseteq Q$ and for
each finite $X \subseteq P$, there exists
$\sigma \in \aut_{A \cup C}(\mathfrak{C})$ such that 
$\sigma(a)/E = g(a/E)$ for each $a \in X$.
This is defined unambiguously since if $x, y \in P$ such 
that $x/E = y/E$ then $\sigma(x)/E = \sigma(y)/E$ for
any automorphism $\sigma \in \aut(\mathfrak{C}/A)$.

We will show first that for any $a, b \models p^n$ and countable
$C \subseteq Q$ there exists $\sigma \in \aut(\mathfrak{C}/A \cup C)$
sending $a$ to $b$.
Next, we will show essentially that the action of $G$ on $P/E$
is $(n+1)$-determined, which we will then use to show that
the action has rank $n$. 
It will follow immediately that $G$ is interpretable in $\mathfrak{C}$,
as in the previous section.
Finally, we will develop the theory of Lascar strong types and strong
automorphisms (over finite sets) to show that $G$ admits hereditarily unique
generics, again, exactly like in the previous section.

We now construct the group more formally.
\begin{definition}
Let $G$ be the group of permutations of $P/E$ such that
for each countable $C \subseteq Q$ and finite $X \subseteq P$
there exists $\sigma \in \aut(\mathfrak{C}/A \cup C)$ such
that $\sigma(a)/E = g(a/E)$ for each $a \in X$.
\end{definition}

$G$ is clearly a group.
We now prove a couple of key lemmas that explain why we chose
$\ecl$ rather than $\bcl$; these will be used to show that
$G$ is not trivial.

\begin{lemma}
Let $a=(a_i)_{i < k}$ be a finite sequence in $P$.
Suppose that $\dim(a/C) = k$, for some $C \subseteq \mathfrak{C}$.
Then there exists $M \prec \mathfrak{C}$ such that 
\[
a_i \not \in \bcl(M a_0 \dots a_{i-1}),
\quad
\text{for each $i < k$}.
\]
\end{lemma}
\begin{proof}
We find models $M_i^j$, for $i \leq j <k$, and
automorphisms $f_j \in \aut(\mathfrak{C}/M_j^j)$ for each $j < k$ such that:
\begin{enumerate}
\item
$A \cup C \cup a_0 \dots a_{i-1} \subseteq M_i^j$ for each $i \leq j < k$.
\item
For each $i < j < n$, $M_i^{j-1} = f_j(M_i^j)$.
\item
$a_j \nonfork_{M_j^j} M_0^j \cup \dots \cup M_{j-1}^j$.
\end{enumerate}

This is possible:  Let $M^0_0 \prec \mathfrak{C}$ 
containing $A \cup C$ be such that
$a_0 \not \in M_0^0$, which exists by definition, and let $f_0$ be the
identity on $\mathfrak{C}$.
Having constructed $M_i^j$ for $i \leq j$, and $f_j$, we let
$M_{j+1}^{j+1} \prec \mathfrak{C}$ contain $A \cup C \cup a_0 \dots a_{j}$
such that $a_{j+1} \not \in M_{j+1}^{j+1}$, which exists by definition.
Let $b_{j+1} \in \mathfrak{C}$ realise $\tp(a_{j+1}/ M_{j+1}^{j+1})$ such that
\[
b_{j+1} \nonfork_{M_{j+1}^{j+1}} M_0^{j} \cup \dots \cup M_j^j.
\]
Such $b_{j+1}$ exists by stationarity of $\tp(a_{j+1}/ M_{j+1}^{j+1})$.
Let $f_{j+1}$ be an automorphism of $\mathfrak{C}$ fixing $M_{j+1}^{j+1})$
sending $b_{j+1}$ to $a_{j+1}$.
Let $M^{j+1}_i = f_{j+1}(M_i^j)$, for $i \leq j$.
These are easily seen to be as required.

This is enough:  Let $M = M_0^{k-1}$.  To see that $M$ is as needed,
we show by induction on $i \leq j < k$, that 
$a_i \not \in \bcl(M_0^j a_0 \dots a_{i-1})$.
For $i = j$, this is clear since 
$a_i \not \in \bcl(M_0^i \cup \dots \cup M_i^i)$.
Now if $j = \ell +1$, $a_i \not \in \bcl(M_0^\ell a_0 \dots a_{i-1})$
by induction hypothesis.  Since $M_0^{\ell+1} = f_{\ell+1}(M_0^\ell)$ and
$f_{\ell+1}$ is the identity on $a_0 \dots a_{i}$, the conclusion follows.
\end{proof}

It follows from the previous lemma that the sequence $(a_i : i < k)$
is a Morley sequence of the quasiminimal type $p_M$,
and hence that (1) it can be extended to any length, and (2) that
any permutation of it extends to an automorphism of $\mathfrak{C}$ over $M$
(hence over $C$).
\begin{lemma} \label{l:ab}
Let $a = (a_i)_{i < n}$ and $b = (b_i)_{i < n}$ be independent
finite sequence in $P$
and a countable $C \subseteq Q$.
Then there exists $\sigma \in \aut(\mathfrak{C}/C)$ such that 
$\sigma(a_i) = b_i$,
for $i < n$.
\end{lemma}
\begin{proof}
By assumption, we have $\dim(a/A \cup C) = \dim(b/A \cup C)$.
By using a third sequence if necessary, we may also assume that
$\dim(ab/A \cup C) = 2n$.
Then, by the previous lemma, there exists $M \prec \mathfrak{C}$
containing $A \cup C$ such that $ab$ is a Morley sequence of $M$.
Thus, the permutation sending $a_i$ to $b_i$ extends to an
automorphism $\sigma$ of $\mathfrak{C}$ fixing $M$ (hence $C$).
\end{proof}

The fact that 
the previous lemma fails for independent sequences of length $n+1$
follows from item (2) of Hypothesis~\ref{h:hyp2}.

We now concentrate on the $n$-action.
We first prove a lemma which is essentially like 
Lemma~\ref{l:1.4}, Lemma~\ref{l:2} and Lemma~\ref{l:3}.
However, since we cannot consider automorphisms fixing all of $Q$,
we need to introduce 
good pairs and good triples.

A pair $(X,C)$ is a {\em good pair} if
$X$ is a countable
infinite-dimensional subset of $P$ with $X = \ecl_A(X) \cap P$;
$C$ is a countable subset of $Q$ such that if
$x_0, \dots, x_n \in X$ with
$x_n \nonfork_A x_0 \dots x_{n-1}$, then there
are $C' \subseteq C$ with 
\[
\dim(x_0 \dots x_n /A \cup C') \leq n,
\] 
$y \in P \setminus \ecl_A(C'x_0 \dots x_{n-1})$ and 
$\sigma \in \aut(\mathfrak{C}/A)$ such that
\[
\sigma(x_n) = y
\quad
\text{ and }
\quad
\sigma(C') \subseteq C.
\]
Good pairs exist; given any countable $X$,
there exists $X' \subseteq P$ countable and $C \subseteq Q$ 
such that $(X', C)$ is a good pair.

A triple $(X,C,C^*)$ is a {\em good triple} if $(X,C)$ is
a good pair, $C^*$ is a countable subset of $Q$ containing $C$,
and whenever two tuples $\bar{a}, \bar{b} \in X$
are automorphic over $A$, then there exists $\sigma \in \aut(\mathfrak{C}/A)$
with $\sigma(\bar{a}) = \bar{b}$ such that, in addition,
\[
\sigma(C) \subseteq C^*.
\]
Again, given a countable $X$, there are $X'$, and $C \subseteq C^*$ such
that $(X',C,C^*)$ is a good triple.

\begin{lemma}
Let $(X,C,C^*)$ be a good triple.
Suppose that $x_0, \dots, x_n \in X$ are independent and
$\sigma(x_i)/E = x_i$ ($i \leq n$) for some $\sigma \in \aut(P/A \cup C^*)$.
Then $\sigma(x)/E = x/E$ for any $x \in X$.
\end{lemma} 
\begin{proof}
We make two claims, which are proved exactly like the stable case
using the definition of good pair or good triple.
We leave the first claim to the reader.
\begin{claim}
Let $(X,C)$ be a good pair. 
Suppose that $x_0, \dots x_{2n-1} \in X$ are independent and 
$\sigma(x_i)/E = x_i/E$ for $i < 2n$ for some 
$\sigma \in \aut(\mathfrak{C}/A \cup C)$.
Then, for all $x \in X \setminus \ecl_A(x_0 \dots x_{2n-1})$ with
$\sigma(x) \in X \setminus \ecl_A(x_0 \dots x_{2n-1})$
we have $\sigma(x)/E = x/E$.
\end{claim}
We can then deduce the next claim:
\begin{claim}
Let $(X,C,C^*)$ be a good triple. 
Suppose that $x_0, \dots, x_n \in X$ are independent and
$\sigma(x_i)/E = x_i/E$ for $i < 2n$ with 
$\sigma \in \aut(\mathfrak{C}/A \cup C^*)$.
Then for each $x \in X \setminus \cl_A(x_0 \dots x_n)$
we have $\sigma(x)/E = x/E$.
\end{claim}
\begin{proof}[Proof of the claim]
Suppose, for a contradiction, that $\sigma(x) \nonfork_A x$.
Using the infinite-dimensionality of $X$ and
the fact that $\sigma(x_i) \in \ecl_A(x_i x_1 \dots x_n)$
we can find $x_i$ for $n \leq i < 2n$
such that
\[
x_i \nonfork_A x x_0 \dots x_{i-1}\sigma(x_1) \dots \sigma(x_{i-1})
\]
and 
\[
\sigma(x_i) \nonfork_A x x_0 \dots x_{i-1}\sigma(x_1) \dots \sigma(x_{i-1}).
\]
It follows that 
\[
x x_0 \nonfork_A x_1 \dots x_{2n} \sigma(x_1) \dots \sigma(x_{2n}),
\]
so there is $\tau \in 
\aut(\mathfrak{C}/A x_1 \dots x_{2n} \sigma(x_1) \dots \sigma(x_{2n}))$
such that $\tau(x) = x_0$.
By definition of good triple, we may assume that $\tau(C) \subseteq C^*$.
Then $\sigma^{-1} \circ \tau^{-1} \circ \sigma \circ \tau$ contradicts
the previous claim.
\end{proof}
The lemma follows from the previous claim by choosing $x_i'$ for $i \leq n$
such that $x_i' \not \in \ecl_A(xx_0 \dots x_n x_0' \dots x_{i-1}')$:
First $\sigma(x_i')/E = x_i'$ for $i \leq n$, and then $\sigma(x)/E = x/E$.
\end{proof}

We now deduce easily the next proposition.

\begin{proposition}\label{p:n+1}
Let $(a_i)_{i \leq n}$ and $(b_i)_{i \leq n}$ be in $P$
such that $\dim((a_i)_{i \leq n}/A) = n+1$.
Let $c \in P$.
There exists a countable $C \subseteq Q$ such that if 
$\sigma, \tau \in \aut(\mathfrak{C}/A \cup C)$ and 
\[
\sigma(a_i)/E = b_i/E = \tau(a_i)/E,
\quad
\text{for each $i \leq n$}
\]
then $\sigma(c)/E = \tau(c)/E$.
\end{proposition}
The value of $\sigma(c)$ in the previous proposition is independent
of $C$.
It follows that the action of $G$ on $P/E$ is $(n+1)$-determined.
We will now show that the action has rank $n$ (so $G$ is automatically
nontrivial).

\begin{proposition}
The action of $G$ on $P/E$ is an $n$-action.
\end{proposition}
\begin{proof}
The $(n+1)$-determinacy of the action of $G$ on $P$ follows from the
previous lemma.  We now have to show that the action has rank $n$.

For this, we first prove the following claim:
If $\bar{a} = (a_i)_{i < n}$ and $\bar{b} = (b_i)_{i < n}$ are in $P$ such
that $\dim(\bar{a}\bar{b}/A) = 2n$ and $c \not \in \ecl(A \bar{a}\bar{b})$, 
then 
there is $d \in P$ such that for each countable
$C \subseteq Q$ there is $\sigma \in \aut(\mathfrak{C}/AC)$
satisfying $\sigma(a_i) = b_i$ (for $i < n$) and $\sigma(c) = d$.

To see this, 
choose $D \subseteq Q$ such that $\dim(\bar{a}c/D) = n$ 
(this is possible by our
hypothesis).  
Suppose, for a contradiction, that no such $d$ exists.
Any automorphism fixing $D$ and sending $\bar{a}$ to $\bar{b}$ must send
$c \in \ecl(A D\bar{b}) \cap P$.
Thus, for each $d \in \ecl(AD\bar{b})$ a countable set $C_d \subseteq Q$
containing $D$
with the property that no automorphism fixing $C_d$ sending $a$ to $b$
also sends $c$ to $d$. 
Since $\ecl(AD\bar{b})$ is countable, we can therefore 
find a countable $C \subseteq Q$ containing $D$ such
that any $\sigma \in \aut(\mathfrak{C}/ A \cup C)$ sending
$\bar{a}$ to $\bar{b}$ is such that $\sigma(c) \not \in \ecl(AD\bar{b})$.
By Lemma~\ref{l:ab}, there does exist $\sigma \in \aut(\mathfrak{C}/ A \cup C)$
such that $\sigma(\bar{a}) = \bar{b}$, and by choice of $D$
we have $\sigma(c) \in \ecl(A D\bar{b})$.
This contradicts the choice of $C$.

We can now show that the action of $G$ on $P/E$ has rank $n$.
Assume that $\bar{a}, \bar{b}$ are independent $n$-tuples of
realisations of $p$.
We must find $g \in G$ such that $g(\bar{a}/E) = \bar{b}/E$.
Let $c \in P \setminus \ecl_A(\bar{a}\bar{b})$ and
choose $d \in P$ as in the previous claim.
We now define the following function $g : P/E \rightarrow P/E$.
For each $e \in P$, choose $C_e$ as in the Proposition~\ref{p:n+1},
{\em i.e.} for any $\sigma, \tau \in \aut(\mathfrak{C}/C_e)$,
such that $\sigma(\bar{a})/E = \bar{b}/E = \tau(\bar{a})/E$
and $\sigma(c)/E = d/E = \tau(c)/E$, we have
$\sigma(e)/E = \tau(e)/E$.
By the previous claim there is $\sigma \in \aut(\mathfrak{C}/C_e)$
sending $ac$ to $bd$.  Let $g(e/E) = \sigma(e)/E$.
The choice of $C_e$ guarantees that this is well-defined.  It is easily 
seen to induce a permutation of $P/E$.
Further, suppose a countable $C \subseteq Q$ is given and a finite
$X \subseteq P$.  Choose $C_e$ as in the previous proposition
for each $e \in X$.  There is $\sigma \in \aut(\mathfrak{C})$ sending
$ac$ to $bd$ fixing each $C_e$ pointwise.
By definition of $g$, we have $\sigma(e)/E = g(e/E)$.
This implies that $g \in G$.
Since this fails for independent $(n+1)$-tuples, by Hypothesis~\ref{h:hyp2},
the action of $G$ on $P$ has rank $n$.
\end{proof}

The next proposition is now proved exactly like Proposition~\ref{c:main}.

\begin{proposition}
The group $G$ is interpretable in $\mathfrak{C}$ (over a finite set).
\end{proposition}

\begin{remark}
Recall that in this case, any complete type over a finite set is equivalent
to a formula (as $\K$ is the class of atomic models of a countable first
order theory).  By $\omega$-homogeneity of $\mathfrak{C}$, for
any finite $B$, any $B$-invariant
is subset of $\mathfrak{C}$ is a countable disjunction of formulas over $A$.
Since the complement of a $B$-invariant set is $B$-invariant,
any $B$-invariant set over a finite set is actually type-definable over $B$.
Hence, the various invariant sets in the above interpretation are all
type-definable over a finite set. 
\end{remark}

It remains to deal with the stationarity of $G$.
As in the previous section, this is done by considering 
{\em strong automorphisms} and {\em Lascar strong types}.
We only need to consider the group of strong automorphisms
over finite sets $C$, which makes the theory easier.

In the excellent case, indiscernibles do not behave as well as 
in the homogeneous case: on the one hand, some indiscernibles cannot
be extended, and on the other hand, it is not clear that a permutation
of the elements induce an automorphism.
However, Morley sequences over models have both of these properties.
Recall that $(a_i : i < \alpha)$ is the {\em Morley sequence} of $\tp(a_0/M)$
if $\tp(a_i/M \{ a_j : j < i \})$ does not split over $M$.
(In the application, we will be interested in Morley sequences inside $P$,
these just coincide with independent sequences.)

We first define Lascar strong types.
\begin{definition}
Let $C$ be a finite subset of $\mathfrak{C}$.
We say that $a$ and $b$ have the same {\em Lascar strong type over $C$},
written $\Lstp(a/C) = \Lstp(b/C)$, if $E(a,b)$ holds for any
$C$-invariant equivalence relation $E$ with a bounded number of classes.
\end{definition}

Equality between Lascar strong types over $C$ is clearly a $C$-invariant
equivalence relation; it is the finest $C$-invariant equivalence relation
with a bounded number of classes.  With this definition, one can 
prove the same properties for Lascar strong types as one has in the
homogeneous case.  The details are in \cite{HyLe:2}; the use of excellence
to extract good indiscernible sequences from large enough sequences is 
a bit different from the homogeneous case, but once one has the fact
below, the details are similar.

\begin{fact}\label{f:main}
Let $I \cup C \subseteq \mathfrak{C}$ be such that $|I|$ is uncountable
and $C$ countable.
Then there is a countable 
$M_0 \prec \mathfrak{C}$ containing $C$ 
and $J \subseteq I$ uncountable
such that
$J$ is a Morley sequence of some stationary type
$p \in S_D(M_0)$.  
\end{fact}

The key consequences are that (1) The Lascar strong types are the
orbits of the group $\Sigma$ of strong automorphisms, and (2) Lascar strong
types are stationary.
We can then show a proposition similar
to Proposition~\ref{p:stat} and Proposition~\ref{p:gen}.
\begin{proposition} $G$ is stationary and 
admits hereditarily unique generics with
respect to $\Sigma$.
\end{proposition}

We have therefore proved:

\begin{theorem}~\label{t:main2}
Let $\K$ be excellent.
Let $\mathfrak{C}$ be a large full model
containing the finite set $A$.
Let $p, q \in S_D(A)$ be unbounded with $p$ quasiminimal.
Assume that there exists an integer $n < \omega$ such that
\begin{enumerate}
\item
For each independent $n$-tuple $a_0, \dots, a_{n-1}$ of realisations
of $p$ and countable $C \subseteq Q$  
we have 
\[
\dim(a_0 \dots a_{n-1} /AC) = n.
\]
\item
For some independent $(n+1)$-tuple $a_0, \dots, a_n$ of realisations
of $p$ and some countable $C \subseteq Q$ we have
\[
\dim(a_0 \dots a_n/AC) \leq n.
\]
\end{enumerate}
Then $\mathfrak{C}$ interprets a group $G$ acting on the geometry $P'$
induced on the realisations of $p$.
Furthermore, either $\mathfrak{C}$ interprets a non-classical group,
or $n \leq 3$ and 
\begin{itemize}
\item
If $n = 1$, then $G$ is abelian and acts
regularly on $P'$;
\item
If $n = 2$, the action of $G$ on $P'$ is isomorphic to the affine
action of $K \rtimes K^*$ on the algebraically closed field $K$.
\item
If $n = 3$, the action of $G$ on $P'$ is isomorphic to the action
of $\PGL_2(K)$ on the projective line $\mathbb{P}^1(K)$ of the
algebraically closed field $K$.
\end{itemize}
\end{theorem}

\begin{question}
Again, as in the stable case, we can produce a group starting from a regular
type only (see \cite{GrHa} for the definition).
Is it possible to get the field ({\em i.e.} hereditarily unique 
generics) starting from a regular, rather than quasiminimal type?
\end{question}


\begin{thebibliography}{99}
\bibitem[Be]{Be} Alex Berenstein, Dependence relations on homogeneous
groups and homogeneous expansions of Hilbert spaces,
Ph.D. thesis, Notre Dame (2002).

\bibitem[Bu]{Bu} Stephen Buechler, {\bf Essential Stability Theory},
Perspective in mathematical logic, Springer-Verlag, Berlin/Heidelberg/New York
(1996).

\bibitem[BuLe]{BuLe} Stephen Buechler and Olivier Lessmann,
{\em Journal of the AMS}, {\bf 6} 1 91--121.

\bibitem[GrHa]{GrHa}  Rami Grossberg  and  Bradd Hart,  
The classification theory of excellent classes, 
   {\em Journal of Symbolic Logic}  {\bf 54} (1989)  1359--1381.

\bibitem[GrLe]{GrLe} Rami Grossberg and Olivier Lessmann,
Shelah's stability spectrum and homogeneity spectrum
in finite diagrams. {\em Archive for mathematical logic},
{\bf 41} (2002) 1, 1--31

\bibitem[Hr1]{Hr} Ehud Hrushovski, Almost orthogonal regular types.
{\em Annals of Pure and Applied Logic} {\bf 45} (1989) 2 139--155.

\bibitem[Hr2]{Hr3} Ehud Hrushovksi, Unidimensional theories are
superstable, {\em Annals of Pure and Applied Logic},
{\bf 50} (1990) 117-138.

\bibitem[Hr3]{Hr1} Ehud Hrushovski,
Stability and its uses, {\em Current development in mathematics},
(1996), Cambridge MA, 61--103.

\bibitem[Hr4]{Hr2} Ehud Hrushovski
Geometric model theory, in {\em Proceedings of the International 
Congress of Mathematicians}, Vol I Berlin (1998).

\bibitem[Hy1]{Hy:4} Tapani Hyttinen, 
On non-structure of elementary submodels of a stable homogeneous structure
{\em Fundament Mathematicae} {\bf 156} (1998) 167--182.

\bibitem[Hy2]{Hy} Tapani Hyttinen,  
Group acting on geometries, in Logic and Algebra, ed. Yi Zhang, 
Contemporary Mathematics,
Vol 302, AMS, 221-233.

\bibitem[Hy3]{Hy:2} Tapani Hyttinen, Finitely generated submodels
of totally categorical homogeneous structures.  Preprint.

\bibitem[Hy4]{Hy:3} Tapani Hyttinen, Characterization for simplicity
in superstable homogeneous structures.  Preprint.

\bibitem[HyLe1]{HyLe} Tapani Hyttinen and Olivier Lessmann,
A rank for the class of elementary
submodels of a superstable homogeneous model.
{\em Journal of Symbolic Logic} {\bf 67} 4 (2002) 1469--1482. 

\bibitem[HyLe2]{HyLe:2} Tapani Hyttinen and Olivier Lessmann,
Independence, simplicity, and canonical bases in excellent classes.
Preprint.

\bibitem[HySh]{HySh} Tapani Hyttinen and Saharon Shelah,
Strong splitting in stable homogeneous models, Annals of Pure and
Applied Logic, {\bf103}, (2000), 201--228.



\bibitem[Ke]{Ke} H. Jerome Keisler, {\bf Model theory for infinitary logic }, 
North-Holland, (1971), Amsterdam.


\bibitem[Le1]{Le:1}  Olivier Lessmann, Ranks and pregometries in 
finite diagrams, {\em Annals of Pure and Applied Logic}, 
{\bf 106}, (2000), pages 49--83.

\bibitem[Le2]{Le:2}  Olivier Lessmann, Homogeneous model theory:
Existence
and Categoricity, in Logic and Algebra, ed. Yi Zhang, 
Contemporary Mathematics,
Vol 302, AMS 149--164.

\bibitem[Le3]{Le:3}  Olivier Lessmann,  Categoricity and U-rank
in excellent classes.  To appear in {\em Journal of Symbolic Logic}.

\bibitem[Ma]{Ma} Angus Macintyre, On $\omega_1$-theories of fields,
{\em Fundamenta Mathematicae}, {\bf 70} (1971) 253--270.
\bibitem[MeSh]{MeSh} Alan Meckler and Saharon Shelah,
$L_{\infty,
\omega}$-free
algebras, {\em Algebra Universalis} {\bf 26} (1989), 351--366.

\bibitem[Pi]{Pi} Anand Pillay, {\bf Geometric Stability Theory}, Oxford
University Press, Oxford (1996).

\bibitem[Po]{Po} Bruno Poizat,  {\bf Groupes Stables}
Nur al-Mantiq wal-Mar'rifah, Villeurbanne (1987).

\bibitem[Re]{Re} Joachim Reineke, Minimale Gruppen, 
{\em Zeitschrift f\"ur Matematische Logik},
{\bf 21} (1975) 357--359.

\bibitem[Sh3]{Sh:3} Saharon Shelah, Finite diagrams stable in power, 
{\em Annals Math. Logic } {\bf 2},
(1970), pages 69--118.



\bibitem[Sh48]{Sh:47} Saharon Shelah, Categoricity in $\aleph_1$ 
of sentences in $L_{\omega_1\omega}(Q)$,
{\em Israel Journal of Math.} {\bf 20} (1975), 127--148.

\bibitem[Sh54]{Sh:54} Saharon Shelah, The lazy model 
theorist's guide to stability, 
{\bf Proc. of a Symp. in Louvain}, March 1975, ed. P. Henrand, 
{\em Logique et Analyse},
18\`eme ann\'ee, {\bf 71-72} (1975), 241--308.


\bibitem[Sh87a]{Sh:87a} Saharon Shelah.
Classification theory for nonelementary classes. I. The number of
uncountable models of $\psi \in L_{\omega_1 \omega}$. Part A.
{\em Israel Journal of Mathematics}, {\bf 46} (1983) 212--240.

\bibitem[Sh87b]{Sh:87b} Saharon Shelah.
Classification theory for nonelementary classes. I. The number of
  uncountable models of $\psi \in L_{\omega _1 \omega}$. Part B.
{\em Israel Journal of Mathematics}, {\bf 46} (1983) 241--273.

\bibitem[Sh]{Sh:a} Saharon Shelah, 
{\bf Classification theory and the number of nonisomorphic models}, 
Rev. Ed., North-Holland, 1990, Amsterdam.

\bibitem[Zi]{Zi1} Boris Zilber, {\bf Uncountably Categorical theories}, 
AMS Translations of mathematical monographs. Vol. 117 (1993).

\bibitem[Zi1]{Zi:2}  Boris Zilber, Covers of the multiplicative group of
an algebraically closed field of characteristic $0$.  Preprint.

\bibitem[Zi2]{Zi:3}  Boris Zilber,  Analytic and pseudo-analytic structures.
Preprint.



\end{thebibliography}
\end{document}